\documentclass[12pt]{amsart}

\usepackage{hyperref}
\usepackage{amsmath}
\usepackage{amsthm}
\usepackage{eucal}
\usepackage{times}
\usepackage{euler}
\usepackage{amssymb}
\usepackage{pst-node}
\usepackage{pst-tree}
\usepackage{verbatim}

\author{Artur Je\.z}
\address{Institute of Computer Science, University of Wroclaw, 
ul.\ Joliot-Curie 15 \\ \mbox{50-384 Wroclaw,} Poland}
\email{Artur.Jez@ii.uni.wroc.pl}

\author{Piotr \'Sniady}
\address{Institute of Mathematics, University of Wroclaw, pl.\ Grunwaldzki 2/4,
\mbox{50-384 Wroclaw,} Poland}

\email{Piotr.Sniady@math.uni.wroc.pl}

\title[Generalized Cauchy identities, trees\dots\ Part II]%
{Generalized Cauchy identities, trees and multidimensional Brownian
motions. \\ Part II: Combinatorial differential calculus}

\theoremstyle{plain}
\newtheorem{lemma}{Lemma}
\newtheorem{analytic-lemma}[lemma]{Analytic lemma}
\newtheorem{theorem}[lemma]{Theorem}
\newtheorem{proposition}[lemma]{Proposition}
\newtheorem{corollary}[lemma]{Corollary}

\newtheorem{inductive}[lemma]{Inductive hypothesis}
\newtheorem{additional-inductive-hypothesis}[lemma]{Additional inductive
hypothesis}

\theoremstyle{definition}
\newtheorem{definition}[lemma]{Definition}

\theoremstyle{remark}
\newtheorem*{remark}{Remark}

\newcommand{\Z}{{\mathbb{Z}}}
\newcommand{\R}{{\mathbb{R}}}
\newcommand{\Lab}{L}
\newcommand{\RR}{\mathcal{R} }
\newcommand{\RRx}{\mathcal{R}[r] }
\newcommand{\E}{\mathbb{E}}
\newcommand{\rut}{r}

\newcommand{\Trees}{{\mathcal{T}}}
 
\newcommand{\homo}{{\E}}
\newcommand{\Graphs}{\mathcal{G}}
\newcommand{\RGraphs}{\mathbb{Z}(\Graphs)}
\newcommand{\derivative}[1]{\frac{\partial}{\partial #1}}
\newcommand{\graf}{\mathbb{G}}
\newcommand{\mniejsze}{<}

\newcommand{\mniejszeD}{\prec}

\DeclareMathOperator{\vol}{vol}

\begin{document}

\begin{abstract}
We present an analogue of the differential calculus in which the role of
polynomials is played by certain ordered sets and trees. Our
combinatorial calculus has all nice features of the usual calculus
and has an advantage that the elements of the considered ordered sets might
carry some additional information. In this way an analytic proof of generalized
Cauchy identities from the previous work of the second author can be directly reformulated in our new
language of the combinatorial calculus; furthermore the additional information
carried by vertices determines uniquely the bijections presented in Part I of
this series.
\end{abstract}

\maketitle


\section{Introduction}

\subsection{Analytic proof of generalized Cauchy identities}
The usual Cauchy identity states that for each nonnegative integer
$k$
\begin{equation}
2^{2k}=\sum_{p+q=k} \binom{2p}{p} \binom{2q}{q},
\label{eq:tozsamosc2}
\end{equation}
where the sum runs over nonnegative integers \begin{math}p,q\end{math}. 
It is a particular case (for $n=2$) of generalized Cauchy identities  
(introduced as a conjecture by 
Dykema and Haagerup \cite{DykemaHaagerup2001}) which state that if $k,n\geq 1$ 
are integers then
\begin{equation}
\label{eq:tozsamoscn} n^{nk}= \text{(number of certain orders on some oriented
trees)},
\end{equation}
where the explicit form of the right-hand side will be recalled later on.
For more on the history of these identities and their applications in theory of
operator algebras 
\cite{Aagaard2003,AagaardHaagerup2004,DykemaHaagerup2001,DykemaHaagerup2003,%
DykemaYan,DykemaJungShlyakhtenko2005} we refer to \cite{Sniady2004BijectionPartI}.

The main idea of the first proof \cite{Sniady2002} of the generalized Cauchy identities
was to associate a
polynomial of a single variable to every oriented tree contributing
to the right-hand side of \eqref{eq:tozsamoscn}. These polynomials for
different values of $n$ turned out to be related by a simple differential
equation and for this reason can be regarded as generalizations of Abel
polynomials. These recursive formulas allowed to express the number of
combinatorial objects contributing to the right-hand side of
\eqref{eq:tozsamoscn} as a certain iterated integral and then to find
explicitly their cardinality.

\subsection{How to convert an analytic proof into a bijection? Combinatorial calculus}

It would be very tempting to find a direct bijection between the orders on the
trees contributing to the right-hand side of \eqref{eq:tozsamoscn} and some
simple set with the cardinality \begin{math} n^{nk} \end{math}. Furthermore, as we pointed out in the
introduction to the Part I of this series \cite{Sniady2004BijectionPartI}, such
a bijection could be used to extract some non-trivial information about
multidimensional Brownian motions and, in particular, to find a
multidimensional analogue of the arc-sine law.

Our idea in looking for such bijective proofs was the following: 
maybe it would be possible to extract the desired
bijection directly from the analytic proof that we had? In the rest of
this article we will present the details of this program.
In general, a map from trees into polynomials of one variable is not 
invertible, i.e.\ usually it is not possible to
extract the original tree from the corresponding polynomial, therefore
we should find an analytic proof of the generalized Cauchy identities which
instead of the differential and integral calculus on polynomials 
uses a differential and integral calculus on richer combinatorial
structures.

In this article we present an analogue of the differential calculus
of one variable in which the role of polynomials is played by
certain ordered sets and trees. In this way an analytic proof of the generalized Cauchy
identities \cite{Sniady2002} can be directly reformulated in our new
language. Furthermore, the additional information carried by the vertices of the graphs
can be used to determine the required bijection uniquely. 

This bijection was already presented in the part I of this series 
\cite{Sniady2004BijectionPartI} in a relatively compact algorithmic way.
A great disadvantage of that approach was that it is by no means clear how the 
bijection from \cite{Sniady2004BijectionPartI} was invented and if it could be generalized
to some other situations. Therefore the main contribution of the current article is not the 
bijection itself but providing a general setup which guides finding such a bijection. 
We also hope that our combinatorial
calculus will be useful in converting  analytic proofs of some other identities into 
bijections.

\subsection{Overview of this article}
This article is organized as follows: in Section \ref{sec:toyexample} we recall the 
polynomials of graphs used in \cite{Sniady2002} to prove its main theorem. In Section \ref{sec:combinatorial-calculus} we refine the ideas of polynomials to much richer combinatorial objects and describe their main properties, which are analogs of the usual laws of calculus for polynomials. 
In Section \ref{sec:problem} we recall the quotient trees, the objects that were studied in
\cite{Sniady2002,Sniady2004BijectionPartI}. 
In Section \ref{sec:towards} we use combinatorial calculus to construct in a conceptual way
the bijection from \cite{Sniady2004BijectionPartI}.
In Section \ref{sec:how-to} we construct explicitly a certain bijection $f$
which plays a crucial role in our construction.

\section{Combinatorial calculus: Toy example}
\label{sec:toyexample}

In this section we present an alternative description of
the usual differential calculus of a single variable. It will serve
us as a toy example after which we shall model our general case.

\subsection{Oriented graphs}

\begin{figure}[td]
\psset{unit=4cm}
\centering{
\begin{pspicture}[](-0.5,-0.1)(1,0.2)
\cnode*(0,0){0.6mm}{pr}
\cnode(0,0){1.2mm}{pr1}
\cnode*(-0.5,0){0.6mm}{pa}
\cnode*(0.5,0){0.6mm}{pb}
\cnode*(1,0){0.6mm}{pc}

\cnodeput*(0,0.1){zz}{$r$}
\cnodeput*(-0.5,0.1){zz}{$a$}
\cnodeput*(0.5,0.1){zz}{$b$}
\cnodeput*(1,0.1){zz}{$c$}

\ncline[doubleline=true]{<-}{pr}{pa}
\ncline[doubleline=true]{<-}{pr}{pb}
\ncline[doubleline=true]{<-}{pb}{pc}

\end{pspicture}}
\caption{Example of an oriented tree. The decorated vertex corresponds to
the root $r$. }
\label{fig:extendingfortree}
\end{figure}
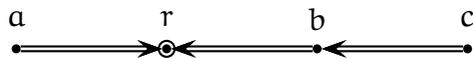

By an \emph{oriented graph} we denote a graph in which every edge is oriented.
We denote an edge from vertex $a$ to a vertex $b$ by $(a,b)$ and write $a \succ
b$ if an edge $(a,b)$ exists. We say that an oriented graph is \emph{acyclic} if
there
is no closed loop of the form $a_1\prec a_2 \prec \cdots \prec a_n \prec a_1$.
In the following all considered oriented graphs will be assumed to be
acyclic, for such graphs the relation \begin{math}\mniejszeD \end{math} can be extended to a
partial order which will be denoted by the same symbol. In other words: $a\mniejszeD
b$ if and only if there is a directed path from $b$ to $a$. An example is presented on Fig.~\ref{fig:extendingfortree}---we can see that $a\succ r$,
$c \succ b\succ r$ and there is no relation between $a$ and $b$ or $a$ and $c$.
In order to avoid ambiguities we will sometimes write $\mniejszeD_G$ in order to make the dependence
on the graph $G$ explicit. 

By a \emph{rooted graph} we denote a graph with a
distinguished vertex, called the \emph{root}, denoted by $r$. The Reader may
restrict attention to the case when the graph $G$ is a
rooted \emph{tree} since this is the case which we consider in this
article.

\subsection{Polynomial associated to a graph}
\label{subsec:polynomial-to-tree}
Let \begin{math} G \end{math} be an oriented graph with a root $r$ and the set of the vertices $V$.
We say that a function \begin{math}f:V\rightarrow [0,1]\end{math} is compatible with the graph $G$ if 
for all pairs of vertices $a,b$ such that \begin{math}a\mniejszeD b\end{math} we also have \begin{math}f(a)\mniejsze f(b)\end{math}.

Let us fix some numbering of non-root vertices. For any fixed \begin{math}x\in[0,1]\end{math} the
set 
\begin{displaymath} Z_G(x)=\big\{ f:V\rightarrow [0,1] \text{ such that $f$ is compatible with
$G$ and $f(r)=x$} \big\}\end{displaymath}
can be identified with a subset of a hypercube \begin{math}[0,1]^{|V|-1}\end{math} and hence its volume
\begin{equation}
\label{eq:definicja-f}
 \E_G(x):= \vol Z_G(x) 
\end{equation}
makes sense. It turns out that \begin{math}\E_G\end{math} is a polynomial of degree \begin{math}|V|-1\end{math}. These
polynomials were the key analytic tool in the proof of generalized Cauchy
identities \cite{Sniady2002}. For technical reasons it is useful to
define $\homo$ also for formal linear combinations of graphs
\begin{displaymath}
\homo_{ \sum_i\! k_i G_i}:=\sum_{i}k_i\ \homo_{G_i},
\end{displaymath}
where $k_i\in\Z$ and $G_i$ is a graph and where $i$ takes a finite number of values.

\subsection{Linear orders}

 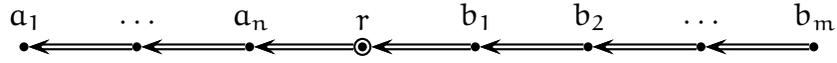
\begin{figure}[td]
 \psset{unit=1.5cm}
 \begin{pspicture}[](-0.5,-0.1)(7.5,.4)
  \cnode*(0,0){0.6mm}{p0}
  \cnode*(1,0){0.6mm}{p1}
  \cnode*(2,0){0.6mm}{p2}
  \cnode*(3,0){0.6mm}{p3}
  \cnode*(4,0){0.6mm}{p4}
  \cnode*(5,0){0.6mm}{p5}
  \cnode*(6,0){0.6mm}{p6}
  \cnode*(7,0){0.6mm}{p7}
  \cnode(3,0){1.2mm}{p3}
  \rput[b](0,0.15){$a_1$}
  \rput[b](1,0.15){\begin{math}\cdots\end{math}}
  \rput[b](2,0.15){$a_n$}
  \rput[b](3,0.15){$r$}
  \rput[b](4,0.15){$b_1$}
  \rput[b](5,0.15){$b_2$}
  \rput[b](6,0.15){\begin{math}\cdots\end{math}}
  \rput[b](7,0.15){$b_m$}
 \ncline[arrowsize=2mm,doubleline=true]{<-}{p0}{p1}
 \ncline[arrowsize=2mm,doubleline=true]{<-}{p1}{p2}
 \ncline[arrowsize=2mm,doubleline=true]{<-}{p2}{p3}
 \ncline[arrowsize=2mm,doubleline=true]{<-}{p3}{p4}
 \ncline[arrowsize=2mm,doubleline=true]{<-}{p4}{p5}
 \ncline[arrowsize=2mm,doubleline=true]{<-}{p5}{p6}
 \ncline[arrowsize=2mm,doubleline=true]{<-}{p6}{p7}
 \end{pspicture}
 \caption{Tree for which the partial order  $ \mniejszeD $ is a total order. It is a $(n,m) $-chain.}
 \label{fig:trywialnygraf}
\end{figure}

In the following we will pay special attention to the case when $G=T$ is a tree for which the
corresponding partial order \begin{math}\mniejszeD\end{math} is a total order. Such a tree must have a
form depicted on Fig.\ \ref{fig:trywialnygraf} for some integers $n,m\geq 0$. We use a name \begin{math}(n,m)\end{math}-chain, or $T_{n,m}$ to denote such a tree.
In this case the set \begin{math}Z_T(x)\end{math} can be viewed as 
\begin{multline}
\label{eq:toyexample} \big\{ (a_1,\dots,a_n,b_1,\dots,b_m): \\ 0<a_1<\cdots<a_n<x<b_1<\cdots<b_m<1\big\} \subseteq [0,1]^{n+m}. 
\end{multline}
We also use a special name for the associated polynomials:
\begin{equation} \E_{n,m}= \frac{x^n}{n!} \frac{(1-x)^m}{m!}. \label{eq:f-nm}\end{equation}
\subsection{Products of graphs}
\label{subsec:product-trees}
\begin{figure}[td]
\psset{unit=2cm}

\begin{pspicture}[](-3,-0.7)(2.5,0.7)
\cnode*(-3,0.5){0.6mm}{pa}
\cnodeput*(-3,0.7){zz}{$a$}
\cnode*(-2.5,0){0.6mm}{pr}
\cnode(-2.5,0){1.2mm}{r}
\cnodeput*(-2.5,0.2){zz}{$r$}
\cnode*(-2,0){0.6mm}{pb}
\cnodeput*(-2,0.2){zz}{$b$}

\cnodeput*(-1.5,0){zz}{\begin{math}\cdot\end{math}}

\cnode*(-1,-0.5){0.6mm}{pe}
\cnodeput*(-1,-0.7){zz}{$e$}
\cnode*(-1,0){0.6mm}{pc}
\cnodeput*(-1,0.2){zz}{$c$}
\cnode*(-0.5,0){0.6mm}{pr2}
\cnode(-0.5,0){1.2mm}{r2}
\cnodeput*(-0.5,0.2){zz}{$r$}
\cnode*(0,-0.5){0.6mm}{pd}
\cnodeput*(0,-0.7){zz}{$d$}

\ncline[doubleline=true]{<-}{r}{pa}
\ncline[doubleline=true]{->}{r}{pb}
\ncline[doubleline=true]{<-}{r2}{pc}
\ncline[doubleline=true]{->}{r2}{pd}
\ncline[doubleline=true]{->}{pc}{pe}

\cnodeput*(0.75,0){zz}{\begin{math}=\end{math}}

\cnode*(1.5,0.5){0.6mm}{pa2}
\cnodeput*(1.5,0.7){zz}{$a$}
\cnode*(1.5,0.){0.6mm}{pc2}
\cnodeput*(1.5,0.2){zz}{$c$}
\cnode*(1.5,-0.5){0.6mm}{pe2}
\cnodeput*(1.5,-0.7){zz}{$e$}
\cnode*(2,0){0.6mm}{pr3}
\cnode(2,0){1.2mm}{r3}
\cnodeput*(2,0.2){zz}{$r$}
\cnode*(2.5,0){0.6mm}{pb2}
\cnodeput*(2.5,0.2){zz}{$b$}
\cnode*(2.5,-0.5){0.6mm}{pd2}
\cnodeput*(2.5,-0.7){zz}{$d$}

\ncline[doubleline=true]{<-}{r3}{pa2}
\ncline[doubleline=true]{->}{r3}{pb2}
\ncline[doubleline=true]{<-}{r3}{pc2}
\ncline[doubleline=true]{->}{r3}{pd2}
\ncline[doubleline=true]{->}{pc2}{pe2}

\end{pspicture}
\caption{Product of trees.}
\label{fig:treemlt}
\end{figure}
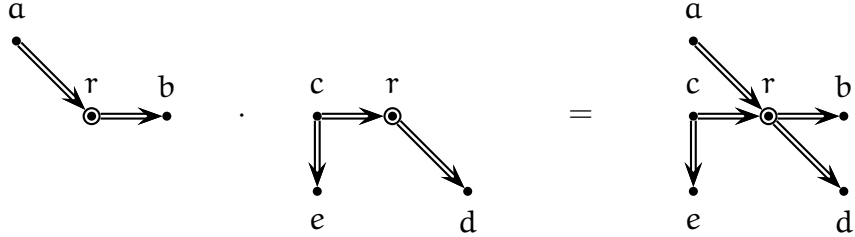

Let $G_1$, $G_2$ be oriented, rooted graphs. For simplicity we shall assume that
the sets of the non-root vertices of these graphs are disjoint. We define a
product $G_1 \cdot G_2$ to be the union $G_1 \cup G_2$ in which the roots of
$G_1$ and $G_2$ are identified. An example of a product of trees is presented on
Fig.~\ref{fig:treemlt}.

We leave it as a simple exercise that 
\begin{equation}
\label{eq:multiplicativity}
 \E_{G_1 \cdot G_2}(x)= \E_{G_1}(x) \ \E_{G_2}(x).
\end{equation}

\subsection{Derivatives of graphs}
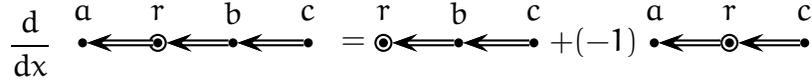
\begin{figure}[td]
\psset{unit=2cm}
\begin{pspicture}[](-3,-0.3)(2.3,0.3)
\cnode*(-2,0){0.6mm}{pr}
\cnode(-2,0){1.2mm}{pr1}
\cnode*(-2.5,0){0.6mm}{pa}
\cnode*(-1.5,0){0.6mm}{pb}
\cnode*(-1,0){0.6mm}{pc}

\cnodeput*(-2,0.2){zz}{$r$}
\cnodeput*(-2.5,0.2){zz}{$a$}
\cnodeput*(-2.85,0){zz}{\begin{math}\displaystyle{\frac{d}{dx}}\end{math}}
\cnodeput*(-1.5,0.2){zz}{$b$}
\cnodeput*(-1,0.2){zz}{$c$}
\cnodeput*(-.7,0){zz}{\begin{math}= \end{math}}

\ncline[doubleline=true]{->}{pr}{pa}
\ncline[doubleline=true]{<-}{pr1}{pb}
\ncline[doubleline=true]{<-}{pb}{pc}
\cnode*(-.5,0){0.6mm}{pr2}
\cnode(-.5,0){1.2mm}{pr12}
\cnode*(0,0){0.6mm}{pb2}
\cnode*(.5,0){0.6mm}{pc2}

\cnodeput*(-.5,0.2){zz}{$r$}
\cnodeput*(0,0.2){zz}{$b$}
\cnodeput*(0.5,0.2){zz}{$c$}

\ncline[doubleline=true]{<-}{pr12}{pb2}
\ncline[doubleline=true]{<-}{pb2}{pc2}

\cnodeput*(0.9,0){zz}{\begin{math}+ (-1) \end{math}}

\cnode*(1.8,0){0.6mm}{pr3}
\cnode(1.8,0){1.2mm}{pr13}
\cnode*(1.3,0){0.6mm}{pa}
\cnode*(2.3,0){0.6mm}{pc}

\cnodeput*(1.8,0.2){zz}{$r$}
\cnodeput*(1.3,0.2){zz}{$a$}
\cnodeput*(2.3,0.2){zz}{$c$}

\ncline[doubleline=true]{->}{pr13}{pa}
\ncline[doubleline=true]{<-}{pr13}{pc}

\end{pspicture}
\caption{Example of derivative of a tree. }
\label{fig:derivative}
\end{figure}

Let $G$ be an oriented rooted graph. We define \begin{math}\frac{d}{dx} G \end{math} to be a formal
linear combination of the graphs given by
\begin{equation}
\label{eq:derivative-first-def}
 \frac{d}{dx} G = \sum_{e} (-1)^e G_{(e)},
\end{equation}
where \begin{math}G_{(e)}\end{math} denotes the graph $G$ with the edge $e$ contracted, \begin{math}(-1)^e=-1\end{math}
if the arrow on $e$ points  towards the root and \begin{math}(-1)^e=1\end{math}
otherwise and the sum runs over all edges $e$ attached to the root. Example is given on Fig. \ref{fig:derivative}.

One can show that
\begin{equation}
\label{eq:rozniczkowanie-ok}
 \E_{\frac{d}{dx} G}(x) = \frac{d}{dx} \E_G(x),
\end{equation}
where the derivative on the right-hand side is the usual derivative of polynomials.

\subsection{Value in $0$ and $1$}
\label{subsec:value-in-0-1}
We focus on the case when $T$ is a chain.
In this case \begin{math}\homo_T=\homo_{n,m}= \frac{x^n}{n!} \frac{(1-x)^m}{m!}\end{math}. The Reader should check that
\begin{equation}\label{simple-embedding} 
\homo_{n,m}(1) = 
\begin{cases} \frac{1}{n!} & \text{if } m=0, \\ 0 & \text{otherwise,} \end{cases}=
\begin{cases}  \sum_{i=0}^{n}\homo_{i,n-i}(x) & \text{if } m=0, \\ 0 & \text{otherwise.} \end{cases}  \end{equation}
The number $\frac{1}{n!}$
in the middle was written as a
linear combination of polynomials on
the right-hand side hence 
we can treat (\ref{simple-embedding}) as a definition of an usual embedding 
$\R\ni a\mapsto a + 0 x + 0 x^2 +\cdots \in \R[x]$. 

We leave it as an exercise to check that analogous result holds true in the case of the value in $0$.

\subsection{Integrals of graphs}
\label{subsec:integrals-first-time}

\begin{figure}[td]
\psset{unit=2.5cm}
\begin{pspicture}[](-2.4,-.65)(2,-.2)
\cnode*(-1.5,-.5){0.6mm}{pr}
\cnode(-1.5,-.5){1.2mm}{pr1}
\cnode*(-2,-0.5){0.6mm}{pa}
\cnode*(-1,-0.5){0.6mm}{pb}
\cnode*(-0.5,-.5){0.6mm}{pc}

\cnodeput*(-1.5,-0.35){zz}{$r$}
\cnodeput*(-2,-0.35){zz}{$a$}
\cnodeput*(-1,-0.35){zz}{$b$}
\cnodeput*(-0.5,-.35){zz}{$c$}

\cnodeput*(0.15,-0.5){zz}{$= (-1) $}
\cnodeput*(-2.3,-0.5){zz}{$\displaystyle{\int_{\#}^{1}}$}
\cnodeput*(-0.3,-0.5){zz}{$dx$}

\cnode*(1,-.25){0.6mm}{pv}
\cnode(1,-.25){1.2mm}{pv1}
\cnode*(1,-0.5){0.6mm}{pr2}
\cnode*(0.5,-.5){0.6mm}{pa2}
\cnode*(1.5,-.5){0.6mm}{pb2}
\cnode*(2,-.5){0.6mm}{pc2}

\cnodeput*(1,-.65){zz}{$x$}
\cnodeput*(1.2,-.2){zz}{$r$}
\cnodeput*(0.5,-.65){zz}{$a$}
\cnodeput*(1.5,-.65){zz}{$b$}
\cnodeput*(2,-.65){zz}{$c$}

\ncline[doubleline=true]{<-}{pr}{pa}
\ncline[doubleline=true]{<-}{pr}{pb}
\ncline[doubleline=true]{->}{pb}{pc}

\ncline[doubleline=true]{<-}{pr2}{pa2}
\ncline[doubleline=true]{<-}{pr2}{pb2}
\ncline[doubleline=true]{->}{pb2}{pc2}
\ncline[doubleline=true]{->}{pr2}{pv1}

\end{pspicture}
\caption{Example of antiderivative. }
\label{fig:graphderivative}
\end{figure}
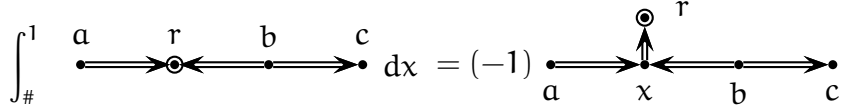

We use the notation that \begin{math} \int_{\#}^{1} f\ dx \end{math} is a function \begin{math}g\end{math} such that
\begin{math} g(y) = \int_{y}^1 f(x)\ dx. \end{math}
We define the corresponding integral for graphs: for a graph $G$ we look for 
$G'$ (which is a formal linear combination of graphs) such that
$\frac{d}{dx}\homo_{G'}=(-1)\cdot \homo_G$ and $\homo_{G'}(1)=0$.
One---particularly elegant---way of constructing such \begin{math}G'\end{math}
is to 
rename the old root of \begin{math} G \end{math} into an ordinary vertex \begin{math}x\end{math} and
to add a new root \begin{math}r\end{math}; then to connect
\begin{math}r\end{math} with \begin{math}x\end{math} with an arrow pointing at
$r$; we denote the resulting graph by $\tilde{G}$, 
see Fig. \ref{fig:graphderivative}. Then
$$\int_{\#}^{1} G\ dx := (-1)\cdot \tilde{G}$$
has the required properties. 

We leave it as an exercise to define the integral $\int_0^{\#} G\ dx$.

\subsection{Extensions of a partial order}
\label{subsec:extensions}
Let $G$ be an oriented graph with the vertex set $V$ and let $\mniejszeD$ be the
corresponding partial order on the set of the vertices. We say that a total
order $\mniejsze$ is compatible with $G$ if $a\mniejszeD b$ implies $a\mniejsze
b$. In this article we are interested in a problem initiated by Dykema and
Haagerup \cite{DykemaHaagerup2001} of studying the set of all total orders
$\mniejsze$ compatible with a given tree $G=T$.

With a small abuse of notation we shall sometimes identify a tree equipped with a total order $(T,\mniejsze)$
with a chain  depicted on Fig.\ \ref{fig:trywialnygraf} with the same vertex
set $V$ and with the order of the vertices $\prec$ specified by $<$. Hence
\begin{equation} \E_{(T,\mniejsze)}= \E_{n,m}= \frac{x^n}{n!} \frac{(1-x)^m}{m!} \end{equation}
makes sense, where $n$ (respectively, $m$) denotes the number of the vertices of $T$ smaller
(respectively, bigger) with respect to $<$ than the root $r$.

One can easily check that
\begin{equation} 
\label{eq:jak-sie-rozklada}
\E_T= \sum_{\mniejsze} \E_{(T,\mniejsze)}, 
\end{equation} 
where the sum runs over all total orders $\mniejsze$ compatible with a tree $T$.

\subsection{Towards the combinatorial calculus}
\label{subsec:towards-combinatorial}

Any polynomial can be written as a linear combination of the polynomials of the form $\E_G$ 
where $G$ is a directed graph. It follows that many operations on polynomials (such as
multiplication, differentiation, integration, taking the value in $0$ or $1$) can be equivalently
performed on the corresponding graphs. Notice, however, that one polynomial can be represented
in many ways as a linear combination of graphs.

Let a directed graph $G$ with a vertex set $V$ be given.
Since the polynomials $\{\E_{n,m}: n+m=|V|-1\}$ defined in
\eqref{eq:f-nm} form a basis of the space of the polynomials of degree at
most $|V|-1$ hence the polynomial $\E_T(x)$ gives us the information about the
number of the total orders $\mniejsze$ compatible with $T$ and such that the number of
the vertices smaller than the root $r$ is specified. Unfortunately, more
detailed information about the order of the vertices with respect to all
possible values of $\mniejsze$ is lost in $\E_T$ and for this reason in Section
\ref{sec:combinatorial-calculus} we shall replace the ring $\R[x]$ of the
polynomials by a richer combinatorial structure.

\section{Combinatorial Differential Calculus}
\label{sec:combinatorial-calculus}

In this section we are going to furnish the objects which appeared in Section
\ref{sec:toyexample} with a richer combinatorial structure. Major changes will
concern the ring of scalars $\R$ in Section \ref{subsec:scalars} and the ring
of polynomials $\R[x]$ in Section \ref{subsec:functions}. In Section 
\ref{subsec:toy-example-revisited} we will revisit the toy example from Section 
\ref{sec:toyexample} and discuss the heuristic meaning of newly introduced
algebraic structures $\RR$ and $\RR[x]$.


\subsection{Graphs}
Let $\Lab$ be a fixed set of labels. 
By $\Graphs$ we denote the set of labeled, directed, acyclic rooted
graphs $G$ such that the root of $G$ is labeled with $\rut\notin\Lab$ 
and all other vertices of $G$ are labeled with different elements of $\Lab$.
By $\Trees$ we denote the subset of $\Graphs$ consisting of directed trees.
We define $\RGraphs$ as the set of formal linear combinations (with integer
coefficients) of elements from $\Graphs$. For graphs in $\Graphs$ we define
multiplication
in the same way as in Section \ref{subsec:product-trees}. By linearity it
extends to $\RGraphs$.

\begin{remark}
In order to avoid technical difficulties we shall always assume that the
vertices of any two graphs involved in any operations considered in the
following have different labels, except for the root, nevertheless this
assumption is not essential.
\end{remark}

\subsection{Scalars}
\label{subsec:scalars}
By $\RR_0$ we denote the set of all finite
sequences with (all different) elements from $\Lab$. We define $\RR$ to be the set of formal linear combinations (with integer
coefficients) of $\RR_0$. This notation was so chosen because $\RR$ is an analogue of the set of
scalars $\R$.

We identify a finite sequence $(a_1,\ldots , a_n)$ with an ordered multiset
$A=\{a_1, \ldots , a_n \} $ where $a_1 \mniejsze \cdots \mniejsze a_n $, also written as $(A,\mniejsze)$.
We can also represent it as a graph with vertices $a_1, \ldots , a_n$ with
oriented edges $(a_n,a_{n-1}), \ldots, (a_2,a_1)$.

Let $A, B \in \RR_0$ be sequences of length $m,n$, respectively. We define $AB\in\RR$ to be the
formal linear combination of $\binom{m+n}{m}$ sequences obtained by
intertwining the sequences $A$ and $B$.
For example:
\begin{multline}(a,b)(c , d )= 
(a,b,c,d) + (a,c,b,d) + \\
(a,c,d,b) + (c,a,b,d) + (c,a,d,b) + (c,d,a,b).
\end{multline}
By linearity this allows us to define the product of two elements of $\RR$.
This multiplication is commutative and associative and it has a unit equal to
the empty sequence $\emptyset$.

\subsection{Polynomials and extensions of partial orders}
\label{subsec:functions}
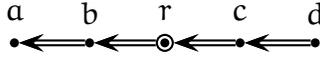
\begin{figure}[td]
\psset{unit=2cm}
\begin{pspicture}[](-1,0)(1,0.3)

\cnode*(-1,0){0.6mm}{a}
\cnodeput*(-1,.2){zz}{$a$}
\cnode*(-0.5,0){0.6mm}{b}
\cnodeput*(-0.5,.2){zz}{$b$}
\cnode*(0,0){0.6mm}{ra}
\cnode(0,0){1.2mm}{r}
\cnodeput*(0,.2){zz}{$r$}
\cnode*(0.5,0){0.6mm}{c}
\cnodeput*(0.5,.2){zz}{$c$}
\cnode*(1,0){0.6mm}{d}
\cnodeput*(1,.2){zz}{$d$}

\ncline[doubleline=true]{<-}{a}{b}
\ncline[doubleline=true]{<-}{b}{r}
\ncline[doubleline=true]{<-}{r}{c}
\ncline[doubleline=true]{<-}{c}{d}
 \end{pspicture}
 \caption{Graphical representation of $(a,b) \otimes (c,d)$.}
\label{fig:porzadek}
\end{figure}

For $A=(a_1, \ldots , a_n ),B=(b_1, \ldots , b_m) \in \RR_0$ we define 
\begin{displaymath}A \otimes B=(a_1, \ldots , a_n , \rut , b_1 , \ldots , b_m )\end{displaymath}
and view this as an ordered multiset (where $\rut$ is some distinguished element
such that $\rut \notin \Lab $), see Fig. \ref{fig:porzadek}.

By $\RRx$ we denote the set of formal linear combinations (with integer
coefficients) of the elements $A \otimes B $ where $A, B \in \RR_0$. This
notation was chosen because $\RRx$ is an analogue of the algebra of polynomials $\R[x]$. 
We replaced the letter $x$ with $r$ in order to stress the connection with the root equipped with label $\rut$.

We equip $\RRx$ with a multiplication by setting 
\begin{displaymath}(A \otimes B)(C \otimes D)
:= (AC) \otimes (BD)\text{.}\end{displaymath}
This product can be described as follows:
for $(A \otimes B)$ and $(C \otimes D)$ we identify the
elements $\rut$ appearing in each of them. The product $(A \otimes B)(C
\otimes D)$ is equal to the formal linear combination of all possible linear
orders on $A \cup B \cup C \cup D \cup \{ \rut \} $ extending the orders on $(A
\otimes B)$ and $(C \otimes D)$ respectively.
Element $\emptyset \otimes \emptyset$ is a unit of this multiplication.
This multiplication is commutative and associative.

Every totally ordered set containing $r$ can be viewed as an element of $\RRx$.
For a directed graph $G\in \Graphs$ with a vertex set $V$ we define $\E_G\in\RRx$ to
be a formal linear combination of all ordered sets $(V,\mniejsze)\in\RRx$, where the
sum runs over all total orders $\mniejsze$ compatible with $G$.

It is easy to check that the following analogue of \eqref{eq:multiplicativity}
holds true.
\begin{proposition}
For all $G_1,G_2\in \RGraphs$
$$\homo_{G_1 \cdot G_2} = \homo_{G_1} \cdot \homo_{G_2},$$ 
where the multiplication on the left-hand side denotes the product of trees and
the multiplication on the right-hand side denotes the product in $\RRx$.
\end{proposition}

\subsection{Derivative}
\label{subsec:derivative}
In analogy to \eqref{eq:derivative-first-def} for $a\in\Lab$ and $G\in
\Graphs$ we define the $a$-derivative by
\begin{equation} 
\label{eq:l-derivative}
\derivative{a} G = G_{(r,a)}-G_{(a,r)}, 
\end{equation}
where $G_{e}$ denotes the graph $G$ in which the edge $e$ was contracted
and non-root label removed or $0$ if edge $e$ does not exist. 
Notice that for simplicity we assumed that
the labels of $G$ are all different therefore there is at most one edge which could be contracted 
in $G_{(r,a)}$, respectively $G_{(a,r)}$; in order to cover the general situation one would have to
consider the formal linear combination of the graphs, each obtained from $G$ by contracting one edge
of the form $(r,a)$, respectively $(a,r)$.

For any $A=(a_1, \ldots , a_n ),B=(b_1, \ldots , b_m) \in \RR_0$ we may treat $A\otimes B\in \RRx$ as an element of $\Graphs$. In this way for any $X \in\RRx$ its derivative 
$\derivative{a} X$ is well-defined and is an element of $\RRx$. This derivative is given explicitly by
\begin{multline*} \derivative{a} \big( A \otimes B \big)=
[a_m=a]\ (a_1,\dots,a_{m-1})\otimes B - \\
 [b_1=a]\  A\otimes (b_2,\dots,b_n).
\end{multline*}

The following result is an analogue of \eqref{eq:rozniczkowanie-ok}.
\begin{proposition}
\label{prop:derlin2}
For every graph $G \in \Graphs$ and $a\in\Lab$
\begin{equation} 
\label{eq:pochodna-jest-ok}
\homo \left(\derivative{a} G\right) = \derivative{a} \homo (G) .
\end{equation}
\end{proposition}
\begin{proof}
This proof is the only place of this paper when we consider graphs which are
not trees or forests. We will use the backward induction with respect to the
number of edges of $G$.

Firstly, let us consider the case when $G$ is a full graph (every pair of
vertices is connected by an oriented edge) and acyclic. Then $\homo(G)$ consists
of exactly one total order on the vertices of $G$ hence there are at most two
summands which contribute to the right-hand side of \eqref{eq:pochodna-jest-ok}.
One can easily check that there are at most two edges $e$ of $G$ adjacent to the
root for which the contracted graph $G_{(e)}$ is acyclic and that they
correspond to the summands on the right-hand side of \eqref{eq:pochodna-jest-ok}
which finishes the proof.

If $G$ is not a full graph we may chose a pair $b,c$ of vertices not connected
by $G$. Let $G_1$ (respectively, $G_2$) denote the graph $G$ augmented by an
edge pointing from $b$ to $c$ (respectively, in the opposite direction). It is
straightforward to show that
$$  \homo \left(\derivative{a} G\right)= \homo \left(\derivative{a} G_1\right)+
\homo \left(\derivative{a} G_2\right) $$
and 
$$ \derivative{a} \homo (G) = \derivative{a} \homo (G_1)+ \derivative{a} \homo
(G_2); $$
notice that it might happen that $G_i$ is not acyclic, in this case
$\homo(G_i)=0$.  The inductive hypothesis can be applied to $G_1$ and $G_2$
which finishes the proof.
%
\end{proof}

The derivative defined above has analogous properties to the usual derivative.
For example, it fulfills Leibniz rule which we leave as a simple exercise.
\begin{proposition}[Leibniz rule] 
\label{prop:Leibnitz}
For every $X_1,X_2 \in \RRx$ and $a \in
\Lab$ we have
$$\derivative{a} (X_1 X_2)=\left( \derivative{a} X_1\right)X_2 +
X_1\left(\derivative{a} X_2\right).$$

For every $ G_1, G_2 \in \Graphs$ and $a \in
\Lab
$ we have
$$\derivative{a} (G_1 G_2)=\left( \derivative{a} G_1\right)G_2 +
G_1\left(\derivative{a} G_2\right).$$
\end{proposition}

\subsection{Embedding of $\RR$ into $\RRx$}
In analogy to embedding (\ref{simple-embedding}) of $\R$ into $\R[x]$ we define embedding $ \iota : \RR \rightarrow \RRx$ given by 
\begin{equation}
\iota (a_1,\ldots,a_n)=\sum_{0\leq k\leq n} (a_1,\ldots,a_k) \otimes
(a_{k+1},\ldots,a_n)
\label{def:iota}
\end{equation}
and which can be alternatively described as follows: to an oriented graph associated to a chain
$(a_1,\ldots,a_n)$ we add an additional vertex $r$, not connected with any other vertices; the 
resulting graph we denote by $G$. Then the right-hand side of \eqref{def:iota} is equal to 
$\E_G$.

The Reader may recognize some similarities of the above definition of $\iota$ to
the definition of embedding \eqref{simple-embedding} of $\homo_{n,m}(1)$ into
$\R[x]$.
 
\begin{proposition}
\label{prop:constant}
For each $X\in\RRx$ the element $\iota(X)$ is a \emph{constant
function} in a sense that for any $a\in \Lab$ 
$$ \derivative{a} \iota(X)=0. $$ 
\end{proposition}

\subsection{Value in $0$ and $1$}
As we have seen in Section \ref{subsec:value-in-0-1} the value in $1$ for chains has 
some nice properties. We define for $A \otimes B \in \RR$
\begin{equation} (A\otimes B) (1)=[ B= \emptyset ]\ \ \iota(A) \in \RRx. 
\label{eq:atone}
\end{equation}
Note that if the root of $G$ has at least one successor then $\big(\E(G)\big)(1)=0$.


Similarly, we define
$$ (A\otimes B) (0)=[ A= \emptyset ]\ \ \iota(B) \in \RRx. $$

\subsection{Integrals}
We have already defined integrals for graphs in  Section
\ref{subsec:integrals-first-time}; we
shall keep this definition with the only change that in the integral of the form
$\int_{\#}^{1} G \ da$ the old root of the graph $G$ will given a label $a$.
We also define an
integral for $\RRx$, namely for $a \in \Lab$ we define
\begin{multline} \int_{\#}^1 (a_1,\ldots,a_m) \otimes (b_1,\ldots,b_n)\ da
= \\ (-1) \sum_{0\leq k\leq m}
(a_1,\ldots,a_k)\otimes(a_{k+1},\ldots,a_m,a,b_1,\ldots,b_n)\in\RRx.
\label{eq:int1}
\end{multline}

The above two integrals (one on $\Graphs$ and one on $\RRx$) are compatible with each other
since 
\begin{multline*} \E\left(\int_{\#}^1 (a_1,\ldots,a_m) \otimes (b_1,\ldots,b_n)\
da\right) 
= \\
\int_{\#}^1 (a_1,\ldots,a_m) \otimes (b_1,\ldots,b_n)\ da, 
\end{multline*}
where on the left-hand side we identify $(a_1,\ldots,a_m) \otimes (b_1,\ldots,b_n)$ with an element of
$\Graphs$ and on the right-hand we treat it as an element of $\RRx$.

Similarly, we define
\begin{multline} \int_0^{\#} (a_1,\ldots,a_m) \otimes (b_1,\ldots,b_n)\ da
= \\ \sum_{0\leq k\leq m}
(a_1,\ldots,a_m,a,b_1,\ldots,b_k)\otimes(b_{k+1},\ldots,b_n)\in\RRx.
\label{eq:int1b}
\end{multline}

\subsection{Basic properties of integrals}

\begin{theorem}[Fundamental theorem of calculus] For any $F \in \RRx$ we have
\label{th:fundamental}
\begin{equation} F = \iota\big[ F(0) \big]+\sum_{a\in \Lab} \int_0^{\#}
\left(\derivative{a} F\right)\ da,
\label{eq:taylor01}
\end{equation}
where the integral should be understood as the integral in $\RRx$.
\end{theorem}
\begin{proof}
Since the integral and the derivative are linear it is sufficient to prove this
equality for $F=A \otimes B = (a_1, \ldots , a_m) \otimes (b_1, \ldots , b_n)$.
Let us calculate the right-hand side in the case when $n,m\geq 1$. We notice,
that although the sum is over all $a \in \Lab$, only $a_m$ and $b_1$ matter:
\begin{multline*} 
\iota\big[F(0)\big]+\sum_{a\in \Lab} \int_0^{\#} \left(\derivative{a} F \right)\
da = \\
\shoveright{0 + \int_0^{\#} (a_1, \ldots ,a_{m-1}) \otimes B\ da_m - 
\int_0^{\#} A \otimes (b_2, \ldots , b_n)\ db_1 = }\\
\shoveleft{\sum_{0 \leq k \leq n} (a_1, \ldots a_{m-1}, a_m, b_1, \ldots , b_k
)\otimes(b_{k+1},\ldots , b_{n}) - } \\
\shoveright{\sum_{0 \leq k \leq n-1} (a_1, \ldots , a_m, b_1, b_2 ,  \ldots ,
b_{k+1}) \otimes (b_{k+2} , \ldots b_n) =} \\
(a_1 , \ldots , a_m) \otimes (b_1 , \ldots , b_n ).
\end{multline*}

The proof in the case when $n$ or $m$ is equal to $0$ requires only minor
modifications.
%
\end{proof}

\begin{corollary}
\label{cor:fundamental}
Let $X,Y\in\RRx$ and $a\in\Lab$. If\/ $Y(0)=0$, $\derivative{a} Y= X$ and 
$\derivative{a'} Y = 0$ holds true for every $a'\neq a $ then 
$$Y= \int_0^{\#} X\ da. $$
\end{corollary}

%


The following proposition shows that the integral is linear with respect to
multiplication by scalars.
\begin{proposition}
For every $C \in \RR$ and $D \in \RRx$ and $ a \in \Lab $ we have
\label{prop:scalmlt}
$$ \iota(C) \cdot \int_0^{\#}D\ da = \int_0^{\#} \iota(C) \cdot D\ da \in\RRx.$$
\end{proposition}
\begin{proof}
We use Corollary \ref{cor:fundamental}
for $Y=\iota(C) \cdot \int_0^{\#}D\ da$ and $X=\iota(C) \cdot D$.
\end{proof}

%

\begin{theorem}[Taylor expansion] 
For any $F \in \RRx$ we have
\begin{equation}
F
= \sum_{\substack{k\geq 0 \\ a_1,\ldots,a_k\in\Lab}}
\iota\left[\left(\derivative{a_1} \cdots \derivative{a_k} F\right)(0) \right] 
\cdot \left( (a_1,\ldots,a_k)\otimes \emptyset \right)
\label{eq:taylor02}
\end{equation}
and
\begin{equation}
F
= \sum_{\substack{k\geq 0 \\ a_1,\ldots,a_k\in\Lab}}
(-1)^k\ 
\iota\left[\left(\derivative{a_1} \cdots \derivative{a_k} F\right)(1) \right] 
\cdot \left( \emptyset \otimes (a_k,a_{k-1},\ldots,a_1) \right).
\label{eq:taylor02prim}
\end{equation}
\end{theorem}
%
%
%
%
%
\begin{proof}
In \eqref{eq:taylor01}
we can further expand $\derivative{a} F$ using the same rule.
But $F$ has a finite number of points and every $\derivative{a} $ reduces the
number of points in $F$, hence the expansion will end after some finite number
of steps and we obtain:
%
\begin{multline*}
F =\sum_{\substack{k\geq 0 \\ a_1,\ldots,a_k\in\Lab}}
 \int_0^{\#} \left[ \cdots \left[\int_0^{\#} \iota\left[ \left( \derivative{a_1}
\cdots
\derivative{a_k} F \right) (0)\right] \ da_1\right] \cdots
\right] da_k = 
\\
\sum_{\substack{k\geq 0 \\ a_1,\ldots,a_k\in\Lab}}
\iota\left[ \left( \derivative{a_1} \cdots \derivative{a_k} F \right) (0)
\right]\ \ \int_0^{\#}
\left[ \cdots \left[\int_0^{\#}  \emptyset \otimes \emptyset\ da_1\right] \cdots
\right] da_k = 
\\
\sum_{\substack{k\geq 0 \\ a_1,\ldots,a_k\in\Lab}} \iota\left[
\left(\derivative{a_1} \cdots \derivative{a_k} F \right)(0) \right] \ \left(
(a_1,\ldots,a_k)\otimes \emptyset\right),
\end{multline*}
where in the second equality we use Proposition \ref{prop:scalmlt}
since $\iota\left[ \derivative{a_1} \cdots \derivative{a_k} F(0)\right]$ is
a scalar and hence we are allowed to move it outside the integral.

The other equation can be proved in an analogous way.
\end{proof}


\subsection{Toy example revisited}
\label{subsec:toy-example-revisited}
It is time to have a look on the definitions introduced in this section and to give them 
heuristic meaning.

As we mentioned in Section \ref{subsec:towards-combinatorial} the polynomial $\E_G\in \R[x]$
fulfills \eqref{eq:jak-sie-rozklada} hence gives some partial information on the 
extensions of the partial order on the vertices of $G$ to total orders.
Similarly, the element $\E_G\in\RRx$  gives (complete) information on such 
extensions; for this reason we regard $\RRx$ as a generalization of $\R[x]$.

Operations $f\mapsto f(0)$ and $f\mapsto f(1)$ map polynomials $\R[x]$ to scalars $\R$;
similarly operations $F\mapsto F(0)$ and $F\mapsto F(1)$ map $\RRx$ to $\RR$, therefore
we regard $\RR$ as an analogue of the set of scalars $\R$.

The usual Taylor expansion for a polynomial $f$ says that
$$ f(x)
= \sum_{k\geq 0}
\left[ \frac{d^k}{dx^k} f(0) \right]  \frac{x^k}{k!}$$ 
and
$$f(x)=\sum_{k\geq 0}
(-1)^k  \left[ \frac{d^k}{dx^k} f(1) \right]  \frac{(1-x)^k}{k!}; $$
clearly each summand on the right-hand sides involves a product of a scalar
$\frac{d^k}{dx^k} f(0) \in \R$ with a polynomial $\frac{x^k}{k!}\in\R[x]$,
respectively of a scalar $(-1)^k \frac{d^k}{dx^k} f(1) \in \R$ with a polynomial
$\frac{(1-x)^k}{k!}\in\R[x]$. In order to define such a product we identify $\R$
with the set of constant polynomials in
$\R[x]$. Similarly, in the Taylor expansion for chains \eqref{eq:taylor02} and
\eqref{eq:taylor02prim} we needed the map $\iota$ to identify $\RR$ as a subset
of $\RRx$. In fact, every derivative $\frac{\partial}{\partial a}$ vanishes on
the image of $\iota$ therefore we can think that $\iota$ maps $\RR$ into
\emph{constant} elements of $\RRx$.

The main difference between calculus on $\R[x]$ and $\RRx$ is the multitude of
derivatives
$\frac{\partial}{\partial a}$ and integrals $\int_0^{\#}  \cdot \ da
$, $\int_{\#}^1  \cdot \ da $ indexed by $a\in L$.
Each of these derivatives and integrals is sensitive to only one label $a\in L$; despite
this multitude we regard our combinatorial calculus as a generalization of the calculus in
\emph{one} variable (as opposite to calculus in several variables) because the total orders 
which we consider have an inherent one-dimensional structure.


\section{Quotient graphs and quotient trees}
\label{sec:problem}
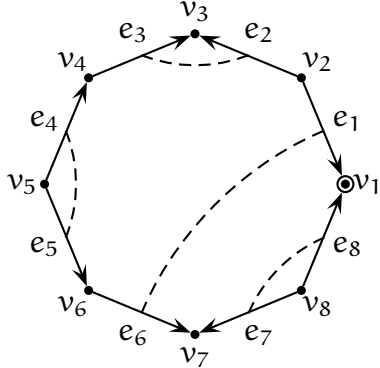
\begin{figure}[td]
\SpecialCoor \psset{unit=0.2cm} \degrees[8]

\begin{pspicture}[](-11,-11)(11,11)
  \cnodeput*(11.5;0){zz}{$v_1$}
 \cnodeput*(11.5;1){zz}{$v_2$}
 \cnodeput*(11.5;2){zz}{$v_3$}
 \cnodeput*(11.5;3){zz}{$v_4$}
 \cnodeput*(11.5;4){zz}{$v_5$}
 \cnodeput*(11.5;5){zz}{$v_6$}
 \cnodeput*(11.5;6){zz}{$v_7$}
 \cnodeput*(11.5;7){zz}{$v_8$}
 \cnodeput*(10.9;0.5){zz}{$e_1$}
 \cnodeput*(10.9;1.5){zz}{$e_2$}
 \cnodeput*(10.9;2.5){zz}{$e_3$}
 \cnodeput*(10.9;3.5){zz}{$e_4$}
 \cnodeput*(10.9;4.5){zz}{$e_5$}
 \cnodeput*(10.9;5.5){zz}{$e_6$}
 \cnodeput*(10.9;6.5){zz}{$e_7$}
 \cnodeput*(10.9;7.5){zz}{$e_8$}
 \cnode*(10;0){0.6mm}{p1}
  \cnode(10;0){1.2mm}{p1}
 \cnode*(10;1){0.6mm}{p2}
 \cnode*(10;2){0.6mm}{p3}
 \cnode*(10;3){0.6mm}{p4}
 \cnode*(10;4){0.6mm}{p5}
 \cnode*(10;5){0.6mm}{p6}
 \cnode*(10;6){0.6mm}{p7}
 \cnode*(10;7){0.6mm}{p8}
 \ncline[arrowsize=2mm]{<-}{p1}{p2}
 \ncline[arrowsize=2mm]{->}{p2}{p3}
 \ncline[arrowsize=2mm]{<-}{p3}{p4}
 \ncline[arrowsize=2mm]{<-}{p4}{p5}
 \ncline[arrowsize=2mm]{->}{p5}{p6}
 \ncline[arrowsize=2mm]{->}{p6}{p7}
 \ncline[arrowsize=2mm]{<-}{p7}{p8}
 \ncline[arrowsize=2mm]{->}{p8}{p1}
 \pnode(9.3;0.5){c1}
 \pnode(9.3;1.5){c2}
 \pnode(9.3;2.5){c3}
 \pnode(9.3;3.5){c4}
 \pnode(9.3;4.5){c5}
 \pnode(9.3;5.5){c6}
 \pnode(9.3;6.5){c7}
 \pnode(9.3;7.5){c8}
 \ncarc[linestyle=dashed,arcangle=-0.5]{c1}{c6}
 \ncarc[linestyle=dashed,arcangle=0.5]{c2}{c3}
 \ncarc[linestyle=dashed,arcangle=0.5]{c4}{c5}
 \ncarc[linestyle=dashed,arcangle=0.5]{c7}{c8}
\end{pspicture}

\caption{Graph $G_{\epsilon}$ corresponding to sequence
$\epsilon=(+1,-1,+1,+1,-1,-1,+1,-1)$. The dashed lines represent a
pairing $\sigma=\big\{ \{1,6\},\{2,3\},\{4,5\},\{7,8\} \} \big\}$.}
\label{fig:graphG}
\end{figure}

The machinery of combinatorial calculus presented in Section
\ref{sec:combinatorial-calculus} was build for the sole purpose of giving
a better understanding of the analytic proof of generalized Cauchy identities
from \cite{Sniady2002}. In this section we introduce the underlying combinatorial structure---quotient graphs and trees.

\subsection{Quotient graphs and quotient trees}
\label{subsec:quetient}


We recall now the construction of Dykema and Haagerup
\cite{DykemaHaagerup2001}. For an integer $k\geq 1$ let $G$ be an
oriented $k$--gon graph with consecutive vertices $v_1,\dots,v_k$
and edges $e_1,\dots,e_k$ (edge $e_i$ connects vertices $v_i$ and
$v_{i+1}$). The vertex $v_1$ is distinguished, see Fig.
\ref{fig:graphG}. We encode the information about the orientations
of the edges in a sequence $\epsilon(1),\dots,\epsilon(k)$ where
$\epsilon(i)=+1$ if the arrow points from $v_{i+1}$ to $v_i$ and
$\epsilon(i)=-1$ if the arrow points from $v_{i}$ to $v_{i+1}$. The
graph $G$ is uniquely determined by the sequence $\epsilon$ and
sometimes we will explicitly state this dependence by using the
notation $G_{\epsilon}$.

Let $\sigma=\big\{ \{i_1,j_1\},\dots,\{i_{k/2},j_{k/2}\} \big\}$ be
a pairing of the set $\{1,\dots,k\}$, i.e.\ pairs $\{i_m,j_m\}$ are
disjoint and their union is equal to $\{1,\dots,k\}$. We say that
$\sigma$ is compatible with $\epsilon$ if
\begin{equation} \epsilon(i)+\epsilon(j)=0 \qquad \text{ for every }\{i,j\}\in\sigma.
\label{eq:orientacja} \end{equation} 
It is a good idea to think that
$\sigma$ is a pairing between the edges of $G$, see Fig.
\ref{fig:graphG}. For each $\{i,j\}\in\sigma$ we identify (or, in
other words, we glue together) the edges $e_i$ and $e_j$ in such a
way that the vertex $v_i$ is identified with $v_{j+1}$ and vertex
$v_{i+1}$ is identified with $v_{j}$ and we denote by $T_{\sigma}$
the resulting quotient graph. The condition \eqref{eq:orientacja}
implies that each edge of $T_{\sigma}$ carries a natural
orientation, inherited from each of the two edges of $G$ it comes
from.

From the following on, we consider only the case when the quotient
graph $T_{\sigma}$ is a tree. One can show \cite{DykemaHaagerup2001}
that the latter holds if and only if the pairing $\sigma$ is
\textit{non--crossing} \cite{Kreweras}; in other words it is not possible
that for some $p<q<r<s$ we have $\{p,r\},\{q,s\}\in\sigma$. The name
of the non--crossing pairings comes from their property that on
their graphical depictions (such as Fig. \ref{fig:graphG}) the
lines do not cross. Let the root $\rut$ of the tree $T_{\sigma}$ be the
vertex corresponding to the distinguished vertex $v_1$ of the graph
$G$.


We say that a sequence $\epsilon=(\epsilon(1),\ldots,\epsilon(n))$ is
\textit{Catalan} when
\begin{equation}
\sum_{i=1}^{k} \epsilon(i) \geq 0\qquad \text { for every $k \leq n$}
\label{def:catalan}
\end{equation}
and 
\begin{displaymath}
\sum_{i=1}^n \epsilon(i) = 0.
\end{displaymath}
We say that $\epsilon$ is \textit{anti-Catalan}, when $-\epsilon$ is Catalan. 
Note, that for a Catalan (respectively: anti-Catalan) sequence and every
non-crossing pairing $\sigma$ in $T_\sigma$ there is no edge incident to the
root and pointing from the root (respectively: towards the root). If
such edge existed then some starting part of $\epsilon$ would sum up to $-1$,
which contradicts (\ref{def:catalan}).

For Catalan sequences one specific pairing will be important in the following
sections, namely the \textit{Catalan pairing}:
\begin{lemma}
For a Catalan (respectively: anti-Catalan) sequence $\epsilon$ there exists a
unique pairing $\sigma$ such that in $T_\sigma$ all edges are directed toward
the root (respectively: in the opposite direction than towards the root). We
call this pairing Catalan pairing.
\label{lem:naturalpairing}
\end{lemma}


\subsection{Preorder}
\label{subsec:Preorder}
\begin{figure}[td]
 \psset{unit=1.65cm}
 \begin{pspicture}[](-2.5,-1)(2.5,0)
 	\pstree[treemode=D,levelsep=0.5,treesep=0.4]{\Tc*{0.6mm} \Tc{1.2mm}~[tnpos=l,tnsep=1mm]{$0$}}{
 		\pstree{\Tc*{0.6mm}~[tnpos=l,tnsep=1mm]{$1$}}{
 			\Tc*{0.6mm}~[tnpos=l,tnsep=1mm]{$2$}
 			\Tc*{0.6mm}~[tnpos=l,tnsep=1mm]{$3$}
 			\Tc*{0.6mm}~[tnpos=l,tnsep=1mm]{$4$}
 			}
 		\Tc*{0.6mm}~[tnpos=b,tnsep=1mm]{$5$}
 		\pstree{\Tc*{0.6mm}~[tnpos=r,tnsep=1mm]{$6$}}{
 			\Tc*{0.6mm}~[tnpos=r,tnsep=1mm]{$7$}
 			\Tc*{0.6mm}~[tnpos=r,tnsep=1mm]{$8$}
 			\Tc*{0.6mm}~[tnpos=r,tnsep=1mm]{$9$}
 			}
 		}
 \end{pspicture}
 \caption{Order $\lhd$ on a rooted plane tree}
 \label{fig:preorder}
\end{figure}
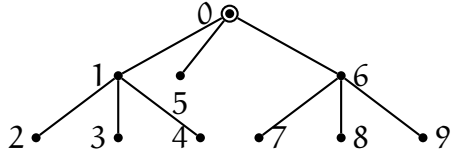

The preorder \cite{StanleyEnumerativeVol2},
denoted by $\lhd$, is defined for trees embedded on
a plane. To obtain it we must traverse
a tree from the root according to the following rule---always choose the
left-most untraversed edge. If there is none, go up the tree. The preorder is
defined by the time of the first visit in the vertex.
Perhaps the illustration on Fig. \ref{fig:preorder} will be better than dwelling this formal definition.

\section{How to convert an analytic proof into a bijection}
\label{sec:towards}


\newcommand{\gwia}{^\ast}
\newcommand{\N}{\mathbb{N}}

\subsection{Formulation of the main result}
\label{subsec:formulation-main-result}
Let 
$T_\epsilon:=\sum_{\sigma} T_\sigma$ be a formal linear combination of the quotient trees, where
$\sigma$ runs over all non-crossing partitions compatible with $\epsilon$.

Let $l_1 \leq l_2 \leq \ldots$ be a weakly increasing sequence of natural numbers.
For each $i\geq 0$ we consider a Catalan sequence
\begin{displaymath}\epsilon_i = \big( \underbrace{(1)}_{l_i \text{times} },
\underbrace{(-1)}_{l_{i-1} \text{times} },
\ldots ,  \underbrace{(-1)^{i-2}}_{l_2 \text{times} },
\underbrace{(-1)^{i-1}}_{l_1 \text{times} }  ,
\underbrace{(-1)^{i}}_{l_1 \text{times} }  ,
\underbrace{(-1)^{i-1}}_{l_2 \text{times} } ,
\ldots , \underbrace{(-1)}_{l_i \text{times} } \big),
\end{displaymath}
where for simplicity instead of $\underbrace{a,a,\dots,a}_{m \text{ times}}$ we
write $\underbrace{a}_{m \text{ times}}$.


The following result was proved (in a slightly different form) in \cite{Sniady2002}.
\begin{theorem}
\label{theo:objetosc}
For each fixed $m \geq 1$ we denote $L=l_1+\cdots+l_m$ and
$\epsilon=\epsilon_m$.
For any $0<x<1$ the value of the polynomial
$\E\big[ T_{\epsilon} \big]\in\R[x]$ is given by
\begin{multline}
\label{eq:objetosc1}
\E\big[ T_{\epsilon} \big](x)=
\vol \Big\{
(x_1,\dots,x_{L})\in\R^{L}: \\
0> x_1 > x_2 > \cdots > x_L > x-m \text{ and }\\
\{x_1,\dots,x_{L},x-m\} \cap [-i,0] \text{ consists of at most } \\
l_1+\cdots+l_i \text{ elements for each } 1\leq i\leq m-1 \Big\}.
\end{multline}

Also
\begin{multline}
\label{eq:objetosc2}
\int_0^1 \E\big[ T_{\epsilon} \big](x) dx=
%
%
%
%
\vol \Big\{
(x_1,\dots,x_{L},z)\in\R^{L+1}: \\ 
0>x_1>x_2>\cdots>x_{L}>z>-m \text{ and }\\
\{x_1,\dots,x_{L},z\} \cap [-i,0] \text{ consists of at most } \\
l_1+\cdots+l_i \text{ elements for each } 1\leq i\leq m-1 \Big\}.
\end{multline}

\end{theorem}


Let $c_1,\dots,c_m$ be different colors.
We consider a function which maps $[-m,0)$ into $[0,1]\times
\{c_1,\dots,c_m\}$
given by  
$$ f(x)= 
\begin{cases} ( \lceil x\rceil-x,c_{-\lfloor x\rfloor})  & \text{if }
(-m-\lfloor x\rfloor) \text{ is odd},
\\
(x-\lfloor x\rfloor, c_{-\lfloor x\rfloor} )  & \text{if } 
(-m-\lfloor x\rfloor) \text{ is even}. \end{cases} $$ 
Notice that the graph of the first coordinate is a zig-zag.

In the following we shall view  $(x,c_i)\in [0,1]\times \{c_1,\dots,c_m\}$ 
as number $x$ decorated with a color $c_i$; in this way the map
$$(x_1,\dots,x_L,z)\mapsto \big\{  f(x_1), \dots, f(x_L), f(z) \big\}$$
provides a bijection between the tuples $(x_1,\dots,x_L,z)$ which contribute to the
set on the right-hand side of \eqref{eq:objetosc2} and sets consisting of $L+1$ elements,
each element being a number from the interval $[0,1]$ and decorated with a color
from the set
$\{c_1,\dots,c_m\}$ with an additional property that for each $1\leq i<m$
at most $l_1+\cdots+l_i$ elements are decorated with colors from the set
$\{c_1,\dots,c_i\}$ (notice that this map is not well-defined or is not a
bijection only on a set
of measure zero with respect to the Lebesgue measure). The element $f(z)$ will
play a special role, we declare it
the root. Therefore Theorem \ref{theo:objetosc} suggests that the following
stronger result could be true.

\begin{theorem}[The main result]
\label{theo:main}
We denote $\epsilon=\epsilon_m$.
There is an explicit way of decorating the vertices of all ordered
trees $(T,<)$ contributing to $\E[T_{\epsilon}]$
with colors $c_1,\dots,c_m$ in such a way that $\E(T_{\epsilon})\in\RRx$ is
equal to
\begin{equation}
\label{eq:definicja-S}
 S=\sum \iota\big[(\underbrace{c_1,\dots,c_1}_{s_1 \text{ times}})\big] \cdots
\iota\big[(\underbrace{c_{m-1},\dots,c_{m-1}}_{s_{m-1} \text{ times}})\big]
\cdot
\big[\emptyset \otimes  (\underbrace{c_{m},\dots,c_{m}}_{s_{m} \text{ times}})  
\big], 
\end{equation}
where the sum runs over tuples of non-negative integers $(s_1,\dots,s_m)$ with a
property that 
$s_1+\cdots+s_m=l_1+\cdots+l_m$
and $s_1+\cdots+s_i\leq l_1+\cdots+l_i$ for $1\leq i\leq m-1$.
\end{theorem}

The remaining part of the article will be devoted to the proof of the above
theorem.

\begin{corollary}
There is an explicit bijection between
\begin{itemize}
\item the set of pairs $(\sigma,\mniejsze)$, where $\sigma$ is a
pairing compatible with $\epsilon$ and $\mniejsze$ is a total order on the
vertices of \/ $T_{\sigma}$ extending $\mniejszeD$;
\item the set of sequences $(a_1,\ldots,a_{L})$ such that
$a_1,\ldots,a_{L}\in\{1,\ldots,m\}$ and for each $1\leq i \leq m-1$ at
most $l_1+\cdots+l_i$ elements of the sequence $(a_n)$ belong to
the set $\{1,\dots,i\}$.
\end{itemize}
\end{corollary}
\begin{proof}
We color each tree $(T,<)$ contributing to $\E[T_{\epsilon}]$ as in Theorem
\ref{theo:main}; additionally we paint the root with the color $c_m$.
To each such a colored tree with a linear order $<$ of the vertices we
associate the sequence of the colors enumerated according to the order $<$;
furthermore each color $c_i$ is replaced by its number $i$. Theorem
\ref{theo:main} shows that it is the required bijection. 
\end{proof}

\subsection{Idea of the proof}
We are going to prove Theorem \ref{theo:main} by induction with respect to $m$.
Let $\epsilon'=-\epsilon_{m-1}$; by reversing the order of the linear orders it
follows that the inductive hypothesis is equivalent to the following statement.

\begin{inductive}
\label{inductive}
There is an explicit way of decorating the vertices of all ordered
trees $(T,<)$ contributing to $\E[T_{\epsilon'}]$
with colors $c_1,\dots,c_{m-1}$ in such a way that $\E(T_{\epsilon'})\in\RRx$ is
equal to
\begin{equation}
\label{eq:definicja-Sprim}
S'=\sum \iota\big[(\underbrace{c_1,\dots,c_1}_{s_1 \text{ times}})\big] \cdots
\iota\big[(\underbrace{c_{m-2},\dots,c_{m-2}}_{s_{m-2} \text{ times}})\big]
\cdot
\big[  (\underbrace{c_{m-1},\dots,c_{m-1}}_{s_{m-1} \text{
times}}) \otimes \emptyset \big],
\end{equation}
where the sum runs over $(s_1,\dots,s_{m-1})$ with a property that
 $s_1+\cdots+s_{m-1}=l_1+\cdots+l_{m-1}$
and $s_1+\cdots+s_i\leq l_1+\cdots+l_i$ for $1\leq i\leq m-2$.
\end{inductive}

A straightforward calculation based on Proposition \ref{prop:Leibnitz} and
Proposition \ref{prop:constant} shows that $S$ given by
\eqref{eq:definicja-S} 
fulfills
\begin{align} 
\label{eq:correct-derivative1}
\left[ \frac{\partial}{\partial x_j} \cdots \frac{\partial}{\partial x_1} S  
\right] (1) & =0
\qquad \text{if } (x_1,\dots,x_j)\neq (c_m,\cdots,c_m), \\
\label{eq:correct-derivative2}
 \Big[ \underbrace{\frac{\partial}{\partial c_m} \cdots 
\frac{\partial}{\partial c_m}}_{j \text{ times}} 
S   \Big] (1) & =\begin{cases} 0 & \text{if } j<l_m,\\
(-1)^{j} \Big[ \underbrace{\frac{\partial}{\partial c_{m-1}} \cdots 
\frac{\partial}{\partial c_{m-1}}}_{j-l_m \text{ times}} 
S' \Big] (1) & \text{if } j\geq l_m.
\end{cases}
\end{align}
Furthermore, the Taylor expansion \eqref{eq:taylor02prim} shows that the
above
derivatives specify $S$ uniquely therefore, in order to prove that some coloring 
is such as claimed by Theorem \ref{theo:main} it is enough to show that  
\eqref{eq:correct-derivative1}--\eqref{eq:correct-derivative2} hold true if
$S$ is replaced by $\E[T_{\epsilon}]$ and---by Inductive hypothesis
\ref{inductive}---if $S'$ is replaced by $\E[T_{\epsilon'}]$, namely
\begin{equation} 
\label{eq:correct-derivative1corrected}
\left[ \frac{\partial}{\partial x_j} \cdots \frac{\partial}{\partial x_1}
\E[T_{\epsilon}]  
\right] (1)  =0
\qquad \text{if } (x_1,\dots,x_j)\neq (c_m,\cdots,c_m), 
\end{equation}
\begin{multline}
\label{eq:correct-derivative2corrected}
 \Big[ \underbrace{\frac{\partial}{\partial c_m} \cdots 
\frac{\partial}{\partial c_m}}_{j \text{ times}} 
\E[T_{\epsilon}]   \Big] (1)  = \\ \begin{cases} 0 & \text{if } j<l_m,\\
(-1)^{j} \Big[ \underbrace{\frac{\partial}{\partial c_{m-1}} \cdots 
\frac{\partial}{\partial c_{m-1}}}_{j-l_m \text{ times}} 
\E[T_{\epsilon'}] \Big] (1) & \text{if } j\geq l_m.
\end{cases}
\end{multline}

This idea of proving Theorem \ref{theo:main} should not come as a surprise since 
basically the same idea appears in the proof of Theorem \ref{theo:objetosc} presented
in \cite{Sniady2002}; namely it was proved there that 
\begin{equation}
\label{eq:stara-rownosc-nie-rdzewieje}
 \Big[ \frac{d^j}{dx^j} \E[T_{\epsilon}]   \Big] (1)  =
\begin{cases} 0 & \text{if } j<l_m,\\
(-1)^{j} \Big[ \frac{d^{j-l_m }}{dx^{j-l_m }} 
\E[T_{\epsilon'}] \Big] (1) & \text{if } j\geq l_m.
\end{cases}
\end{equation}
with the only difference that here $\E[T_{\epsilon}]\in\R[x]$ is just a polynomial
and the derivatives are the usual derivatives.

In the following we are going to analyze the analytic proof of 
\eqref{eq:stara-rownosc-nie-rdzewieje}
presented in \cite{Sniady2002} and find its ramifications in our more general context.

\subsection{Orders of derivatives}
\label{subsec:derivatives}

The following lemma was critical in the proof of \eqref{eq:stara-rownosc-nie-rdzewieje}
in \cite{Sniady2002}.

\begin{analytic-lemma}
\label{analytic-lemma}
Let $\epsilon=\big(\epsilon_1,\dots,\epsilon_{2k} \big)$ be a Catalan sequence.
Then
\begin{displaymath} 
\frac{d^k}{dx^k} \E_{T_{\epsilon}} = (-1)^k, \end{displaymath}
where $\E_{T_{\epsilon}}\in\R[x]$ is the usual polynomial in one variable as
defined in \eqref{eq:definicja-f}. In the light of
\eqref{eq:rozniczkowanie-ok} it can be equivalently stated that
\begin{equation} 
\label{eq:pochodna-catalana}\frac{d^k}{dx^k} {T_{\epsilon}} = (-1)^k \cdot
(\text{trivial tree consisting only of the root}),
\end{equation}
where the derivative of trees should be understood as in
\eqref{eq:derivative-first-def}. 
\end{analytic-lemma}
We are going to find an analogue of the combinatorial identity
\eqref{eq:pochodna-catalana} in which the robust derivative
\eqref{eq:derivative-first-def} would be replaced by a more refined derivative
\eqref{eq:l-derivative}. We will do it by labeling the vertices in such a way
that the derivatives of the form 
$\frac{\partial}{\partial x_k} \cdots  \frac{\partial}{\partial x_1} T_{\epsilon}$ 
would have a particularly simple structure for any $x_1,\dots,x_l\in L$. 
This labeling will turn out to be the one required by Theorem \ref{theo:main}.

Applying $\frac{d^k}{dx^k}$ to a tree is a sequential removing of $k$
edges in all possible ways. In order to keep track of all such ways we use
the notion of the \emph{order of derivatives}.

\begin{definition}
$\Delta=(x_1, \ldots , x_n )$ is an \textit{order of the derivatives} of a tree $T$,
when $x_1, \ldots , x_n$ are different vertices of $T$ (all different from the root) and
$x_k$ is a neighbor of the root in $\derivative{x_{k-1}} \ldots \derivative{x_1}
T$. In this case we write $\partial_{\Delta}T=\derivative{x_n} \ldots
\derivative{x_1} T $.

\label{def:naturaldifferencing}
\end{definition}
This name \emph{order of derivatives} is motivated by the fact that 
\begin{displaymath} \left(\sum_{x \in L}\derivative{x}\right)^n T= 
\sum_{\Delta=(x_1, \ldots , x_n )}  \partial_\Delta T \end{displaymath}
and so we may view $\Delta$ as an order of removing vertices from the tree.

Since all the vertices have different labels we sometimes refer to $\Delta$ as
tuple of labels. We are mainly interested in the case when $\Delta$ consists of
all vertices different from the root; in this case we call it total
order of derivatives. 
Also we will write $x_i \mniejsze_{\Delta} x_j$ if $i<j$;
in other words the derivative
$\frac{\partial}{\partial x_i}$ is applied
in $\partial_\Delta$ before the derivative $\frac{\partial}{\partial x_j}$.

\subsection{Orders of derivatives and involutions on trees}

A careful analysis (which can be found in Section \ref{sec:how-to}) 
shows that the proof of the Analytic lemma \ref{analytic-lemma} 
presented in  \cite{Sniady2002} is based on finding cancellations between all summands contributing
to \eqref{eq:pochodna-catalana}. In fact, these cancellations arise from some implicit pairing 
between all possible orders of derivatives.
The latter statement is formalized in the following lemma.  

\begin{lemma}
\label{lem:pairingofderivatives}
Let a Catalan sequence $\epsilon$ be fixed. Let us consider pairs $(\sigma,\Delta )$ 
where $\sigma$ is a non-crossing pairing compatible with $\epsilon$ and 
$\Delta$ is a total order of derivatives  on $T_\sigma$. 
If we exclude the case when $\sigma$ is the Catalan pairing and $\Delta$ is the preorder, 
then there is an explicit involution $f$ without fix-points 
on the rest of those pairs such that
\begin{equation} \big[f(\sigma_1,\Delta_1)=(\sigma_2,\Delta_2)\big] \text{ implies }
   \big[\partial_{\Delta_1} T_{\sigma_1} = -
\partial_{\Delta_2} T_{\sigma_2}\big].
\label{eq:siblings}
\end{equation}

Also, if $\Delta_1=(x_1^{(1)}, \ldots , x_n^{(1)})$ is such that for some $k$
the vertices $x_1^{(1)}, \ldots, x_k^{(1)}$ satisfy
\begin{itemize}
\item $x_1^{(1)} \lhd  x_2^{(1)} \lhd \cdots \lhd x_k^{(1)}$,
\item $ x_1^{(i)} \succ r$ for $i=1 , \ldots, k$
\end{itemize}
then the same applies to $\Delta_2=(x_1^{(2)}, \ldots , x_n^{(2)})$, namely
\begin{itemize}
\item $x_1^{(2)} \lhd  x_2^{(2)} \lhd \cdots \lhd x_k^{(2)}$,
\item $ x_i^{(2)} \succ r$ for $i=1 , \ldots, k$.
\end{itemize} 
\end{lemma}
The explicit form of the pairing $f$ will be constructed in Section
\ref{sec:how-to} where we will also prove that it fulfills the above properties.


Any pair $(T_\sigma,\Delta)$, where $\Delta$ is a total order of
derivatives on $T_\sigma$ will be called a configuration. 

\subsection{Trees with orders of derivatives versus trees with total orders}

We encounter here a major difficulty, namely the labeling wanted in Theorem 
\ref{theo:main} is a labeling of the vertices of trees $T_{\sigma}$ equipped with
compatible total orders while the analytic proof of \cite{Sniady2002} and Lemma 
\ref{lem:pairingofderivatives} suggest that we should rather work with configurations
$(T_\sigma,\Delta)$. As one can easily see, these two notions are quite different and 
in the following we will have to reconcile them.


In order to do this let us have a look on the following problem. 
Let an oriented tree $T$
be fixed.
\emph{Can we associate some canonical order of derivatives
$\Delta=(x_1,\dots,x_k)$ 
on $T$ to a given  total order $<$ on the vertices of\/ $T$?} Let
$a_1,\dots,a_m\prec r$ 
and $b_1,\dots,b_n\succ r$ be the direct neighbors of the root $\rut$; 
we may assume that their numbering was chosen in such a way that
$a_1<\cdots<a_m<\rut<b_1<\cdots<b_n$. 
In the extreme case when $k=1$ and 
the tree $T$ has no other vertices than
$a_1,\dots,a_m,\rut,b_1,\dots,b_n$
we may think that the total order $<$ gives to the set of the vertices a structure of a 
chain, cf Fig.\  \ref{fig:trywialnygraf}. In this chain the root $\rut$ has (at most) two
direct neighbors which could be differentiated, namely $a_m$ and $b_1$. 
Therefore it seems reasonable to say that an order of derivatives $\Delta=(x_1)$ is compatible
with the total order $<$ if $x_1\in\{a_m,b_1\}$. The above discussion motivates heuristically
the following definition.
\begin{definition}
\label{def:compatible}
Let $T$ be an oriented tree, $<$ be a total order on the vertices of $T$ which is compatible
with the orientations of the edges and let $\Delta=(x_1,\dots,x_k)$ be an order of derivatives.
We say that $<$ and $\Delta$ are compatible if for every $1\leq l\leq k$ one of the 
following conditions holds true:
\begin{itemize}
\item $x_l$ is the biggest element (with respect to $<$) of the set
$ \big\{ x\in T_{l}: x \mniejszeD_{T_{l}} r \big\} $,
\item $x_l$ is the smallest element (with respect to $<$) of the set
$ \big\{ x\in T_{l}: r \mniejszeD_{T_{l}} x \big\} $,
\end{itemize}
where the tree $T_l$ is given by
$T_l=\frac{\partial}{\partial x_{l-1}} \cdots \frac{\partial}{\partial x_{1}} T $.
\end{definition}

We are going to investigate which information is preserved 
when we replace a total linear order by one of the corresponding total orders of
derivatives. The answer to this problem will be given in Theorem
\ref{theo:induces1} below.



\begin{definition}
Let $x\neq r$ be a vertex of an oriented  tree $T\in\Trees$ and let
$(x_0,\dots,x_n)$ be the shortest path connecting the root $r=x_0$ and $x=x_n$. 
Let $1\leq i\leq n-1$ be the biggest index such that $x_{i-1} \prec x_i \succ
x_{i+1}$ or $x_{i-1} \succ x_i \prec x_{i+1}$. If such an index exists we say
that $x_i$ is \emph{the last bend} of $x$. 
We say that vertices $x,y \in T$ are in the same layer if one of the following conditions holds true:
\begin{itemize}
\item $x$ and $y$ have the same last bend; 
\item $x,y\prec r$ or $x,y\succ r$.
\end{itemize}

\end{definition}
\begin{definition}
Let a tree $T$ be given. Let $\Delta_1$ and $\Delta_2$ be total orders of
derivatives. 
We say that $\Delta_1 \sim_d \Delta_2$ if $x \mniejsze_{\Delta_1} y
\iff x \mniejsze_{\Delta_2} y$ holds for all $x,y$ which are in the
same layer. 

Let $\mniejsze_1,\mniejsze_2$ be total linear orders on the vertices of $T$
compatible with the orientations of the edges.
We say that $\mniejsze_1\ \sim_o\ \mniejsze_2$ if $x \mniejsze_1 y
\iff x \mniejsze_2 y$ holds for all $x,y$ which are in the same
layer.
\end{definition}



\begin{theorem}
\label{theo:induces1}
The notion of compatibility (Definition \ref{def:compatible}) provides a
bijection between the classes of the equivalence relation $\sim_o$ and the 
classes of the equivalence relation $\sim_d$. 
\end{theorem}
\begin{proof} The proof is immediate. \end{proof}

Let some equivalence class $c$ of the equivalence relation $\sim_o$ 
be given (hence it is an equivalence class of the relation $\sim_d$); we define
$$ \E[T,c]= \sum_{< \in c} T_< \in \RRx,$$
where $T_<\in\RRx$ denotes the chain defined by the linear order $<$. The same
definition makes sense if $c$ is replaced by an equivalence class of $\sim_d$.

\begin{theorem}
\label{theo:induces2}
Let $c$ be an equivalence class of\/ $\sim_d$ and $\Delta=(x_1,\dots,x_m)$ be an
order of derivatives
of the tree $T$. If\/ $\Delta$ is not a prefix of any element of\/ $c$ then
$$\partial_\Delta \E[T,c]= 0. $$
If\/ $\Delta$ is a prefix of a total order of derivatives $(x_1,\dots,x_n)\in c$
then
$$\partial_\Delta \E[T,c]= 
\E\bigg(\partial_\Delta T, \big[(x_{m+1},\dots,x_n)\big]_{\sim_d} 
\bigg).$$
\end{theorem}
\begin{proof}
Assume for simplicity that $\Delta$ is a total order of derivatives.
We consider any linear order $<$
which belongs to the class $c$. We consider a graph $T'$ with
the same set of the vertices as $T$ and for each pair of vertices $v<w$ we
draw an arrow from $w$ to $v$ if one of the following conditions
hold true:
\begin{itemize}
\item vertices $v$ and $w$ are in the same layer and 
there is no vertex $z$ in the same layer which fulfills $v<z<w$;
\item one of the vertices (let us denote it by $p$) is the last bend of the
other (let us denote it by $q$) and there is no vertex $z$ in the same layer as
$q$ which fulfills $v<z<w$. 
\end{itemize}
It is easy to check that the above definition does not depend on the choice of
$<$ and that the resulting graph $T'$ is a tree. Furthermore, each total linear
order $<$ is compatible with the orientations of the edges of $T'$  if and only
if $<$ belongs to $c$. It follows that
$$\E[T,c]=\E[T'].$$
Furthermore, any total order of derivatives $\Delta$ belongs to the equivalence
class specified by $c$ if and only if $\Delta$ is an order of derivatives on
$T'$.
Therefore
$$\partial_\Delta \E[T,c]=\E[\partial_\Delta T']$$
which finishes the proof.

The case when $\Delta$ is not necessarily total follows in a similar way.
\end{proof}

\subsection{Towards the labeling}
\label{subsec:towards-the-labeling}

Our ultimate goal is to find some special labeling (coloring) of the vertices of 
the trees $(T_{\sigma},<)$ equipped with total orderings; in order to do this we
shall follow the following two principles:
\begin{enumerate}
\item \label{cond:we_require_nice1}
for pairs $(T_{\sigma},<)$ from the same equivalence class of $\sim_o$
the labeling of the vertices of $T_{\sigma}$ should be the same 
(our motivation is Theorem \ref{theo:induces2} since 
this requirement would imply that for any configuration
$(T_{\sigma},\Delta)$ the coloring of the vertices of $T_{\sigma}$ is
well-defined);

\item \label{cond:we_require_nice2}
for any two configurations $\left(T_{\sigma_i},(x_1^{(i)},\dots,x_n^{(i)}) 
\right)$, $i\in\{1,2\}$, which are paired by Lemma
\ref{lem:pairingofderivatives} we require that their colorings should be
compatible in a sense that
the color of $x_k^{(1)}$ in tree $T_{\sigma_1}$ should coincide with the color
of $x_k^{(2)}$ in tree $T_{\sigma_2}$ for any value of $k\in\{1,\dots,n\}$.
\end{enumerate}


%

As we shall see in the following, these conditions are quite restrictive and
there are not too many such labelings.


Let us consider a graph $\graf$ with the set of vertices equal to the set of
pairs $(T_{\sigma},[\Delta]_{\sim_d})$ or, equivalently, the set of pairs
$(T_{\sigma}, [<]_{\sim_o} )$. If 
$$f(T_{\sigma},\Delta) =
(T_{\sigma'},\Delta')\quad \text{and}\quad \big|\{v \in T_{\sigma}: v \succ
r\}\big| > \big|\{v \in T_{\sigma'}: v \succ r\}\big|$$ 
then 
we put a directed edge $\big( (T_{\sigma},[\Delta]_\sim),(T_{\sigma'},[\Delta
']_\sim) \big)$ between the appropriate equivalence classes. We denote the
connected component of the vertex $(T_{\sigma},[\Delta]_{\sim_d})\in\graf$ by
$\graf(T_{\sigma},[\Delta]_{\sim_d})$ or simply $\graf(T_{\sigma},\Delta)$. 

\begin{proposition}
\label{prop:sinks_exist}  
Assume that $\epsilon$ is a Catalan sequence. Every vertex of\/ $\graf$ has at
most one outgoing edge and $\graf$ is acyclic therefore every component has a
\emph{sink}---the only vertex with no outgoing edges. 

Configuration $(T_{\sigma},\Delta)$ with $\Delta=(x_1,\dots,x_n)$ corresponds
to a sink $(T_{\sigma},[\Delta]_{\sim_d})$  if
and only if there exists a number $k$ with a property that
$\{x_1,\dots,x_k\}=\{v\in T_{\sigma}: v \succ r \}$ and the order of the
elements $(x_1,\dots,x_k)$ coincides with the preorder.
\end{proposition}
We postpone the proof to Section \ref{sec:how-to}.

\begin{theorem} 
\label{theo:sink_is_enough1}
We assume that $\epsilon$ is a Catalan sequence.
Let $(T_{\sigma},[\Delta]_{\sim_d})$ be a sink and let some coloring of the
vertices of\/
$T_{\sigma}$ be given (respectively, let some coloring of the vertices of all
sinks be given). Then there is
a unique way of extending this coloring to the vertices of the trees
$(T_{\sigma'},\Delta')$ (or, equivalently, trees $(T_{\sigma'},<')$) which
contribute to $\graf(T_{\sigma},\Delta)$ (respectively, to all such trees) so
that conditions
(\ref{cond:we_require_nice1}), (\ref{cond:we_require_nice2}) from the
beginning of Section \ref{subsec:towards-the-labeling} are fulfilled.
\end{theorem}
\begin{proof}
Proposition \ref{prop:sinks_exist} shows that there is a unique path in the
graph $\graf$ which connects any vertex with the sink. On the other hand,
requirement (\ref{cond:we_require_nice2}) from the beginning of Section
\ref{subsec:towards-the-labeling} shows that if the vertices of $\graf$
corresponding to the trees $(T_{\sigma_1},\Delta_1)$, $(T_{\sigma_2},\Delta_2)$
are connected by an edge then the coloring of vertices of
$(T_{\sigma_1},\Delta_1)$ uniquely determines the coloring of the vertices of
$(T_{\sigma_2},\Delta_2)$ which finishes the proof.
\end{proof}

Let a coloring of the vertices of all trees $(T_{\sigma},<)$ be given as above.
We consider the element
\begin{equation}
\label{eq:definicja_grupowania}
\E[ T_{\epsilon} ] = \sum_{(T_\sigma,<)} (T_{\sigma},<)
\in \RRx, 
\end{equation}
where each summand $(T_{\sigma},<)$ is identified with the corresponding chain
of the colors of the vertices of $T_{\sigma}$.
The following theorem shows that thanks to this coloring the derivatives of $\E[
T_{\epsilon} ]$ have a particularly simple structure.

\begin{theorem}
\label{theo:sink_is_enough2} Assume that $\epsilon$ is a Catalan sequence. 
Then for any colors $d_1,\dots,d_k$
\begin{equation} 
\label{eq:cudne-kasowanie}
\frac{\partial }{\partial d_k} \cdots \frac{\partial }{\partial d_1}
\E[ T_{\epsilon} ] = 
\sum_{(T_{\sigma},[\Delta]_{\sim_d})} \E \left[ \frac{\partial }{\partial x_k}
\cdots \frac{\partial }{\partial x_1} (T_{\sigma},\Delta) \right],
\end{equation}
where the symbol $(T_{\sigma},\Delta)$ on the
right-hand side denotes the tree $T_{\sigma}$
with the coloring of the vertices specified by $(T_{\sigma},\Delta)$ as in
Theorem \ref{theo:sink_is_enough2};
the sum runs over equivalence classes $(T_{\sigma},[\Delta]_{\sim d})$ for
which there exists a representative $\Delta=(x_1,\dots,x_n)$
with the following properties:
\begin{itemize} 
\item the colors of the vertices $(x_1,\dots,x_k)$
of the tree $(T_{\sigma},\Delta)$ are equal to $(d_1,\dots,d_k)$;
\item all edges removed from
$T_{\sigma}$ by $\frac{\partial }{\partial x_k} \cdots \frac{\partial
}{\partial x_1}$ are oriented towards the root;
\item all edges removed from
$T_{\sigma}$ by $\frac{\partial }{\partial x_k} \cdots \frac{\partial
}{\partial x_1}$ are removed in the order which coincides with the preorder. 
\end{itemize}


\end{theorem}

\begin{proof}
We may group the summands on the right-hand side of
\eqref{eq:definicja_grupowania} according to the equivalence classes of
$\sim_d$ and apply Theorem \ref{theo:induces2} therefore
$$\frac{\partial }{\partial d_k} \cdots \frac{\partial }{\partial d_1}
\E[ T_{\epsilon} ] = 
\sum_{(T_{\sigma},[\Delta]_{\sim_d})} \E \left[ \frac{\partial }{\partial x_k}
\cdots \frac{\partial }{\partial x_1} (T_{\sigma},\Delta) \right],$$
where the sum on the right-hand side runs over all
$(T_{\sigma},[\Delta]_{\sim_d})$, $\Delta=(x_1,\dots,x_n)$ such that the colors
of the vertices $(x_1,\dots,x_k)$
of the tree $(T_{\sigma},\Delta)$ are equal to $(d_1,\dots,d_k)$.

Let us disregard for a moment all vertices of the trees other than
$x_1,\dots,x_k$; we denote such a truncated tree by $\tilde{T}$. The involution
$f$ from Lemma \ref{lem:pairingofderivatives} can be applied to the set of such
truncated trees; the contribution of all paired trees cancel and the only
remaining trees are as prescribed in the formulation of the theorem.

For completeness of this proof one should check that that the involution $f$
applied to a the truncated tree $(\tilde{T},(x_1,\dots,x_k))$ gives the
truncation of $f(T,\Delta)$. This, however, will become obvious in Section 
\ref{sec:how-to} when the explicit form of $f$  will be given.
\end{proof}

\subsection{Proof of the main result}

In fact, we shall prove Theorem \ref{theo:main} together with the following
additional inductive hypothesis.

\begin{additional-inductive-hypothesis}
\label{additional-hypothesis}
If $(T_{\sigma},\Delta)$ is a configuration with $\Delta=(x_1,\dots,x_n)$ and
$p$ is such that all edges removed by the derivative $\frac{\partial}{\partial
x_p} \cdots \frac{\partial}{\partial x_1} T_{\sigma}$ are oriented towards the
root and they are removed in the order which coincides with the preorder then
the vertices $x_1,\dots,x_p$ are painted by color $c_m$.
\end{additional-inductive-hypothesis}

\begin{proof}[Proof of Theorem \ref{theo:main}]
As we already pointed out in Section \ref{subsec:formulation-main-result}, it
is enough to find a coloring with a property that
equations \eqref{eq:correct-derivative1corrected},
\eqref{eq:correct-derivative2corrected} are fulfilled.

In view of Theorem \ref{theo:sink_is_enough1} it is enough to define the
coloring on the sinks of the graph $\graf$. Equation \eqref{eq:cudne-kasowanie}
shows that
\begin{equation}
\label{eq:gdynia}
\left(\frac{\partial }{\partial d_k} \cdots \frac{\partial
}{\partial d_1}
\E[ T_{\epsilon} ]\right)(1) =
\sum_{(T_{\sigma},[\Delta]_{\sim_d})} \E \left[ \frac{\partial }{\partial x_k}
\cdots \frac{\partial }{\partial x_1} (T_{\sigma},\Delta) \right](1),
\end{equation}
where the sum runs over equivalence classes $(T_{\sigma},[\Delta]_{\sim d})$ for
which there exists a representative $\Delta=(x_1,\dots,x_n)$ with a property
that $\{x_1,\dots,x_k\}=\{v\in T_{\sigma}: v \succ r \}$ and which fulfills the
conditions from Theorem \ref{theo:sink_is_enough2}. In other words:
$(T_{\sigma},[\Delta]_{\sim_d})$ contributes to the above sum if and 
only if $(T_{\sigma},[\Delta]_{\sim_d})$ corresponds to a sink in $\graf$ for
which $k$ is as prescribed in Proposition \ref{prop:sinks_exist}.


Let us fix some sink $(T_{\sigma},[\Delta]_{\sim_d})$; we can always assume 
that $\Delta=(x_1,\dots,x_n)$ is as specified in Proposition
\ref{prop:sinks_exist} and that $k$ is as above. Each of the edges of
$T_{\sigma}$ arises from a pair of
the edges of the polygonal graph $G_{\epsilon}$; let us unglue the edges which
are removed by the
derivative $\frac{\partial }{\partial x_k} \cdots \frac{\partial }{\partial
x_1}$ and let us decorate these edges. The resulting graph is a polygon with
some trees attached or, in other words, it is the quotient graph $T_{\sigma'}$
where the partial pairing $\sigma'$ is subset of the pairing $\sigma$. Among
unglued edges there must be the edges of $G_{\epsilon}$ which correspond to the
$l_i$ initial and the $l_i$ final elements of $\epsilon$; let us remove these
$2l_i$ edges. The resulting graph $T_{\sigma''}$ is a quotient graph of the
polygonal graph $G_{\epsilon'}$. We can view $T_{\sigma''}$ as some polygonal
graph $G_{\delta}$ with some trees attached. We glue the edges of
$G_{\delta}$ by the Catalan pairing, we denote the resulting tree by
$T_{\sigma'''}$. This tree has $k-l_i$ decorated edges which we denote by
$(y_1,\dots,y_{k-l_i})$, in the order given by the preorder; we denote by
$\Delta'=(y_1,\dots,y_{k-l_i},x_{k+1},x_{k+2},\dots,x_n)$ an order of
derivatives on $T_{\sigma'''}$. 

Now we are ready to define the coloring of the vertices of the tree
$T_{\sigma}$: we paint the vertices $x_1,\dots,x_k$ with color $c_m$; all
other
vertices we paint with the same colors as the corresponding vertices in the
tree $(T_{\sigma'''},\Delta')$ (this coloring is given by the inductive
hypothesis).




In order to show that this coloring indeed fulfills
\eqref{eq:correct-derivative1corrected} we use 
Eq.\ \eqref{eq:gdynia}.
Since by the construction of the coloring all vertices
in the set  $\{x: x \succ r\}$ of the sink are colored by color $c_m$ therefore
there is no summand which contributes to the right-hand side of
\eqref{eq:gdynia} which finishes the proof of
\eqref{eq:correct-derivative1corrected}.

We will prove now that \eqref{eq:correct-derivative2corrected} holds true.
In order to get a non-zero value in $1$ of the right-hand side of
\eqref{eq:gdynia} one has to remove by derivatives all vertices $\{x: x \succ
r\}$. Since there are at least $l_m$ such vertices
therefore the first part of \eqref{eq:correct-derivative2corrected} follows.

For the second part of \eqref{eq:correct-derivative2corrected}
note that by Additional inductive hypothesis \ref{additional-hypothesis}
vertices $y_1,\dots,y_{k-l_i}$ have color $c_{m-1}$ therefore
$(T_{\sigma'''},\Delta')$ is one of the summands which contribute to \eqref{eq:gdynia}
 applied to the derivative
$\underbrace{\frac{\partial}{\partial c_{m-1}} \cdots 
\frac{\partial}{\partial c_{m-1}}}_{j-l_m \text{ times}}\E[T_{\epsilon'}]$.
The corresponding relation is one to one, that is given a summand contributing to the above sum we can find a sink such that
it is its $(T_{\sigma'''},\Delta')$. 

It remains now to prove that the coloring costructed above fulfills Additional
inductive hypothesis \ref{additional-hypothesis}. In order to do this assume
that in $(T,\Delta)$ all edges removed by $\frac{\partial}{\partial
x_p} \cdots \frac{\partial}{\partial x_1} T_{\sigma}$ are oriented towards the
root and removed in the order coinciding with the preorder. By Lemma
\ref{lem:pairingofderivatives} the same happens in the corresponding sink
$(T',\Delta')$ and by the construction of the coloring all vertices in the set
$\{ v : v \succ r \}$ removed in the sink are painted by color $c_m$ which
finishes the proof.

In this way we proved that the presented coloring has the required properties.
One can check that this coloring coincides with the one presented in the paper
\cite{Sniady2004BijectionPartI}. 
\end{proof}

\section{How to convert an analytic proof into involution $f$}
\label{sec:how-to}

In this section we will analyze the proof of Analytic lemma \ref{analytic-lemma}
presented in \cite{Sniady2002}  and we will show how the involution $f$ can be
constructed out of it.

\subsection{Sketch of proof of Analytic lemma \ref{analytic-lemma}}
\label{subsec:sketch}

The proof of Analytic lemma \ref{analytic-lemma} presented in
\cite{Sniady2002} 
was based on the following observation. 
The left-hand side of 
\eqref{eq:pochodna-catalana} is a sum over all possible ways of choosing a
quotient
tree $T_{\sigma}$ and then choosing the order of the derivatives
$\frac{\partial}{\partial x_k} \cdots  \frac{\partial}{\partial x_1}$.
Let us concentrate on the last derivative $\frac{\partial}{\partial x_k}$ in the
above
product. This derivative removes the edge $x_k$ of the tree $T_{\sigma}$; 
this edge of $T_{\sigma}$ corresponds to a pair of edges in the polygonal graph
$G_{\epsilon}$. Since one of the ends of $x_k$ is a leaf therefore
the corresponding pair of edges $e_i,e_{i+1}$ must be adjacent and 
have opposite orientations, i.e.\ $\epsilon_i+\epsilon_{i+1}=0$. 
We denote by $(\epsilon)_i$ the sequence 
$\epsilon$ with elements $\epsilon_i,\epsilon_{i+1}$
corresponding to these edges removed. It is easy to see that the 
contribution to \eqref{eq:pochodna-catalana} of all summands for which the 
edge $x_k$ is fixed is equal to $-\epsilon_i \frac{d^{k-1}}{dx^{k-1}}
T_{(\epsilon)_i}$. 
In this way we proved that
\begin{equation} 
\label{eq:smart-sniady}
\frac{d^k}{dx^k} T_{\epsilon} = 
\sum_{\substack{i: \\ \epsilon_i+\epsilon_{i+1}=0}}
(-1) \epsilon_i\cdot  \frac{d^{k-1}}{dx^{k-1}} T_{(\epsilon)_i}.
\end{equation}

The assumption that $\epsilon$ is a Catalan sequence implies that 
each sequence $(\epsilon)_i$ is Catalan as well and the inductive hypothesis can
be applied;
it follows that
$$ \frac{d^k}{dx^k} T_{\epsilon} = 
\sum_{\substack{i: \\ \epsilon_i+\epsilon_{i+1}=0}}
(-1) \epsilon_i \cdot (-1)^{k-1}.$$
Now it is enough to notice that if $i_1<\dots<i_l$ are all indices such that
$(\epsilon_{i_j},\epsilon_{i_j+1})=(1,-1)$ and $i'_1<\dots<i'_{l'}$
are all indices such that $(\epsilon_{i'_j},\epsilon_{i'_j+1})=(-1,1)$ then the
sequences
$(i_j)$ and $(i'_j)$ are interlacing: $i_1<i'_1<i_2<\cdots<i'_{l-1}<i_l$
therefore
$$ \frac{d^k}{dx^k} T_{\epsilon} = 
\sum_{\substack{i: \\ \epsilon_i+\epsilon_{i+1}=0}}
\epsilon_i \cdot (-1)^k= \big[ l \cdot (+1) + (l-1) \cdot (-1) \big]
(-1)^k=(-1)^k, $$
which finishes the proof.

\subsection{How to find the involution in Lemma \ref{lem:pairingofderivatives}}
\label{subsec:how-to-find-involution}

We are going to find the involution as in Lemma \ref{lem:pairingofderivatives}
by following the proof presented in Section \ref{subsec:sketch}. 

Firstly, since the proof in Section \ref{subsec:sketch} is inductive,
our construction of the involution $f$ will be inductive as well.
For the shortest possible Catalan sequence $\epsilon=(1,-1)$ there is 
only one pairing $\sigma$ (which is the Catalan pairing)
and only one total order of derivatives $\Delta$ (which coincides with preorder
$\lhd$)
therefore there is nothing to pair and $f=\emptyset$. 

If $\epsilon$ consists of at least four elements, the proof of
\eqref{eq:smart-sniady}
suggests that we should group all pairs 
\begin{equation}
\label{eq:march2007}
 \big(T_\sigma, \Delta= (x_1,\dots,x_k) \big)
\end{equation}
into classes according to the value of the edge $x_k$ or, in other words,
according to the pair of edges $(e_i,e_{i+1})$.
Similarly as in Section \ref{subsec:sketch} the set of pairs
\eqref{eq:march2007} 
for a fixed value of $i$ can be identified with the set of pairs 
$\big(\sigma',
(x_1,\dots,x_{k-1}) \big)
$
which contribute to the the derivative $\frac{d^{k-1}}{dx^{k-1}}
T_{(\epsilon)_i}$.
It follows that for each group of pairs $(T_\sigma,\Delta)$ corresponding to a
given
value of $i$ (such that $\epsilon_i+\epsilon_{i+1}=0$)
we can construct a pairing $f$ inductively. Nevertheless, for each value of $i$
there is one pair $(T_\sigma,\Delta)$ which remains unpaired and in order
to finish the construction of $f$ we should find some pairing between these
remaining elements.
The requirement \eqref{eq:siblings} implies that that an index $i$ for which
$(\epsilon_i,\epsilon_{i+1})=(1,-1)$ should be paired with an index $i'$
for which $(\epsilon_i,\epsilon_{i+1})=(-1,1)$.

As we already mentioned in Section \ref{subsec:sketch},  
if $i_1<\dots<i_l$ are all indices such that
$(\epsilon_{i_j},\epsilon_{i_j+1})=(1,-1)$ and $i'_1<\dots<i'_{l-1}$
are all indices such that $(\epsilon_{i'_j},\epsilon_{i'_j+1})=(-1,1)$ then the
sequences
$(i_j)$ and $(i'_j)$ are interlacing: $i_1<i'_1<i_2<\cdots<i'_{l-1}<i_l$.
We have a relative freedom in choosing the pairing between the elements of the
set
$\{i_1,\dots,i_{l}\}$ and $\{i'_1,\dots,i'_{l-1}\}$ and we decided to choose 
a fairly natural one: $i_1\leftrightarrow i'_1, \dots,  i_{l-1}\leftrightarrow
i'_{l-1}$
with the index $i_l$ unpaired.

The above inductive procedure determines the involution $f$ claimed in Lemma
\ref{lem:pairingofderivatives} uniquely, nevertheless this description is quite
implicit
and we will present its explicit form in the following section.

\subsection{Explicit form of the involution $f$}
\label{subsec:explicit_form}

\begin{lemma}
\label{lem:jak-wyglada-sierota}
The pair $(T_\sigma,\Delta)$ for which $\sigma$ is the Catalan pairing and
$\Delta$ is
the preorder $\lhd$ is the unique pair which is unpaired by the involution $f$
described
in Section \ref{subsec:how-to-find-involution}.
\end{lemma}
\begin{proof}
The element \eqref{eq:march2007} unpaired by $f$ must belong to the class
considered in Section \ref{subsec:how-to-find-involution} corresponding to
the only unpaired index $i_{l}$ for which $(\epsilon_{i_l},
\epsilon_{i_l+1})=(1,-1)$
therefore the edge $x_k$ must be oriented towards the root.
Also the pair
\begin{equation}
\label{eq:short-product}
 \big(T_\sigma',
(x_1,\dots,x_{k-1}) \big)
\end{equation}
corresponding to the sequence $(\epsilon)_{i_l}$ must be unpaired by $f$;
by iterating this argument it follows that all edges in the tree $T_{\sigma}$
are
oriented towards the root hence $\sigma$ is the Catalan pairing.

Since the only element unpaired by $f$ must belong to the class
considered in Section \ref{subsec:how-to-find-involution} corresponding to
the only unpaired index $i_{l}$ therefore $x_k$ is the last edge with respect to
the
preorder $\lhd$ in the tree $T_{\sigma}$. By passing to \eqref{eq:short-product}
and iterating the argument we see that $\Delta$ coincides with $\lhd$.
\end{proof}

Let a pair $(T_\sigma, \Delta)$ be given, we shall compute explicitly 
the value of $f(\sigma, \Delta)$. It is easy to see that the algorithm described
implicitly in Section \ref{subsec:how-to-find-involution} looks for the smallest
value $p$ for which the pairing $f$ is well-defined for
$(\sigma',\frac{\partial}{\partial x_{p}} \cdots \frac{\partial}{\partial
x_1})$ 
(where $\sigma'$ is a pairing for the sequence $\epsilon$ with the elements
corresponding to the edges $x_{p+1},x_{p+2},\dots$ removed). Lemma 
\ref{lem:jak-wyglada-sierota} shows that $p$ is the smallest value for which 
one of the following conditions hold:
\begin{enumerate}
\item 
\label{case:bad-case-1}
the order of the vertices in the sequence $(x_1,\dots,x_p)$ does not
coincide with the preorder; 
\item 
\label{case:bad-case-2}
edge removed by the last derivative in the product
$\frac{\partial}{\partial x_p} \cdots \frac{\partial}{\partial x_1}$ is not
oriented towards the root.
\end{enumerate}

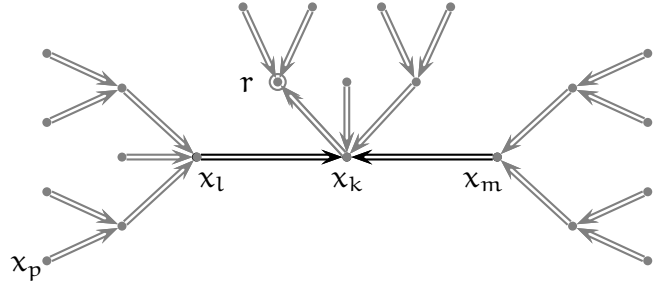
\begin{figure}
 \psset{unit=2cm}
 \begin{pspicture}[](-1,-0.8)(1,1.2)
 \cnode*(0,0){0.6mm}{A}
 \cnode*(1,0){0.6mm}{C}
 \cnode*(-1,0){0.6mm}{B}
 \ncline[arrowsize=2mm,doubleline=true]{<-}{A}{B}
 \ncline[arrowsize=2mm,doubleline=true]{<-}{A}{C}
 \rput[t](0,-0.1){$x_k$}
 \rput[t](0.9,-0.1){$x_m$}
 \rput[t](-0.9,-0.1){$x_l$}

\psset{linecolor=gray}

\renewcommand{\psedge}[2]{\ncline[doubleline=true,arrowsize=2mm]{->}{#1}{#2}}
 \rput[B](0,0){\pstree[treemode=U,levelsep=0.5,treesep=0.4]{\Tc*{0.6mm}}{
     \pstree{ \Tc{1.2mm}~[tnpos=l,tnsep=2mm]{$r$} \Tc*{0.6mm}  }{
    
\renewcommand{\psedge}[2]{\ncline[doubleline=true,arrowsize=2mm]{<-}{#1}{#2}}
   \Tc*{0.6mm}  \Tc*{0.6mm}}
   \renewcommand{\psedge}[2]{\ncline[doubleline=true,arrowsize=2mm]{<-}{#1}{#2}}
   \Tc*{0.6mm}\pstree{\Tc*{0.6mm}}{\Tc*{0.6mm} \Tc*{0.6mm}}}}

 \renewcommand{\psedge}[2]{\ncline[doubleline=true,arrowsize=2mm]{<-}{#1}{#2}}
 \rput[B](-1.6,0){\pstree[treemode=L,levelsep=0.5,treesep=0.4]{\Tc*{0.6mm}}{
     \pstree{\Tc*{0.6mm}}{\Tc*{0.6mm} \Tc*{0.6mm}}
\Tc*{0.6mm}\pstree{\Tc*{0.6mm}}{\Tc*{0.6mm}
     \Tc*{0.6mm}~[tnpos=l,tnsep=0mm]{$x_p$}    }}}

 \rput[B](1.5,0){\pstree[treemode=R,levelsep=0.5,treesep=0.4]{\Tc*{0.6mm}}{
     \pstree{\Tc*{0.6mm}}{\Tc*{0.6mm}
     \Tc*{0.6mm}} \pstree{\Tc*{0.6mm}}{\Tc*{0.6mm} \Tc*{0.6mm} }}}

 \end{pspicture}
\caption{Tree for which the order of the vertices given by $\Delta$
does not coincide with the preorder. The indices fulfill $k<l<m<p$. Only the
vertices $r,x_1,\dots,x_p$ were shown.}

 \label{fig:transform0}
\end{figure}

We shall concentrate in the following on the case (\ref{case:bad-case-1}).
Let us traverse the plane tree $T_{\sigma}$ truncated to the vertices
$r,x_1,\dots,x_p$ starting from the vertex $x_p$. We denote by $(x_m,x_k)$ the
first edge we traverse in the direction opposite to its orientation and by
$(x_l,x_k)$ the previous edge, cf Figure \ref{fig:transform0}.
 
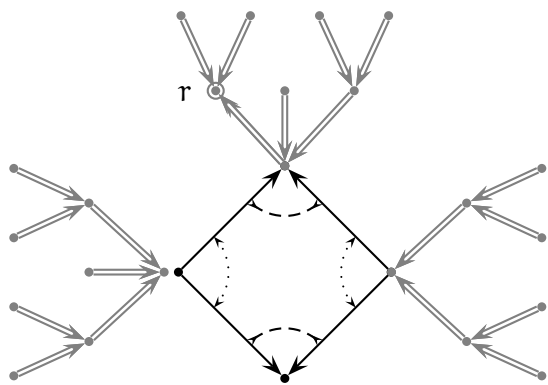
\begin{figure}
\psset{unit=2cm}
\begin{pspicture}[](-1,-0.8)(1,1.8)
 \cnode*(0.7071067,0){0.6mm}{E}
 \cnode*(0,0.7071067){0.6mm}{N}
 \cnode*(-0.7071067,0){0.6mm}{W}
 \cnode*(0,-0.7071067){0.6mm}{S}
 \psarc[linestyle=dashed]{>-<}(0,0.7071067){0.33333}{225}{315}
 \psarc[linestyle=dashed]{>-<}(0,-0.7071067){0.33333}{45}{135}
 \psarc[linestyle=dotted]{<->}(0.7071067,0){0.33333}{135}{225}
 \psarc[linestyle=dotted]{<->}(-0.7071067,0){0.33333}{-45}{45}
 \ncline[arrowsize=2mm,doubleline=false]{<-}{N}{E}
 \ncline[arrowsize=2mm,doubleline=false]{<-}{N}{W}
 \ncline[arrowsize=2mm,doubleline=false]{<-}{S}{E}
 \ncline[arrowsize=2mm,doubleline=false]{<-}{S}{W}

\psset{linecolor=gray}

\renewcommand{\psedge}[2]{\ncline[doubleline=true,arrowsize=2mm]{->}{#1}{#2}}
 \rput[B](0,0.7071067){\pstree[treemode=U,levelsep=0.5,treesep=0.4]{\Tc*{0.6mm}}
{
     \pstree{ \Tc{1.2mm}~[tnpos=l,tnsep=2mm]{$r$} \Tc*{0.6mm}  }{
    
\renewcommand{\psedge}[2]{\ncline[doubleline=true,arrowsize=2mm]{<-}{#1}{#2}}
   \Tc*{0.6mm}  \Tc*{0.6mm}}
   \renewcommand{\psedge}[2]{\ncline[doubleline=true,arrowsize=2mm]{<-}{#1}{#2}}
   \Tc*{0.6mm}\pstree{\Tc*{0.6mm}}{\Tc*{0.6mm} \Tc*{0.6mm}}}}

 \renewcommand{\psedge}[2]{\ncline[doubleline=true,arrowsize=2mm]{<-}{#1}{#2}}
 \rput[B](-1.303,0){\pstree[treemode=L,levelsep=0.5,treesep=0.4]{\Tc*{0.6mm}}{
     \pstree{\Tc*{0.6mm}}{\Tc*{0.6mm} \Tc*{0.6mm}}
\Tc*{0.6mm}\pstree{\Tc*{0.6mm}}{\Tc*{0.6mm}
     \Tc*{0.6mm}~[tnpos=l,tnsep=0mm]{
}    }}}

 \rput[B](+1.2071067,0){\pstree[treemode=R,levelsep=0.5,treesep=0.4]{\Tc*{0.6mm}
}{
     \pstree{\Tc*{0.6mm}}{\Tc*{0.6mm}
     \Tc*{0.6mm}} \pstree{\Tc*{0.6mm}}{\Tc*{0.6mm} \Tc*{0.6mm}}}}
 \end{pspicture}
 \caption{The tree from Figure \ref{fig:transform0} after ungluing the edges
$(x_l,x_k)$, $(x_m,x_k)$.}
 \label{fig:transform1}
\end{figure}

The edge of $T_{\sigma}$ removed by the last derivative
$\frac{\partial}{\partial x_p}$ in the product
$\frac{\partial}{\partial x_p} \cdots \frac{\partial}{\partial x_1}$
corresponds to a pair of edges $(\epsilon_{i_j},\epsilon_{i_j+1})$ of the
polygonal graph $G_{\epsilon}$. These two edges were paired with 
$(\epsilon_{i'_j},\epsilon_{i'_j+1})$ which correspond to one of the half-edges
constituting $(x_l,x_k)$ and to one of the half edges constituting $(x_m,x_k)$.
The pairing $f$ is defined in the following, seemingly complicated way: we
unglue all edges appearing in the tree $T_{\sigma}$ truncated to the vertices
$r,x_1,\dots,x_p$, we glue together one of the half-edges $(x_l,x_k)$ to one of
the half edges $(x_m,x_k)$ and then we glue all remaining unglued edges by the
Catalan pairing. Notice, however, that this operation can be equivalently
described in a much simpler way: we unglue two edges of $T_{\sigma}$, namely 
$(x_l,x_k)$ and $(x_m,x_k)$, cf Figure \ref{fig:transform1} and we reglue them
in a different way, cf Figure \ref{fig:transform2}. The resulting tree
$T_{\sigma'}$ is the wanted tree such that $f(T_{\sigma},\Delta) =
(T_{\sigma'},\Delta')$. 

In order to describe the order of derivatives $\Delta'$
it will be convenient to label the vertices of $T_{\sigma'}$ with the same
labels as the vertices of $T_{\sigma}$, namely $x_1,x_2,\dots$. Question arises
therefore: how to label the vertices of $T_{\sigma'}$ in such a way that
$\Delta'=(x_1,x_2,\dots)$. There is a seemingly complicated way of describing
this labeling: firstly, we need to bother only with the part
of the tree $T_{\sigma}$ which consists of the vertices
$r,x_1,\dots,x_p$ and the corresponding part 
of the tree $T_{\sigma'}$ (we do not change the gluings of any other edges and
the labels given to other vertices remain the same in $T_{\sigma'}$). 
Secondly, in the process of calculating the derivative
$\frac{\partial}{\partial x_p} \cdots \frac{\partial}{\partial x_1}
T_{\sigma'}$ the first $p-1$ derivatives should remove the edges oriented
towards the root in the order given by the preorder and the derivative
$\frac{\partial}{\partial x_p}$ should remove the only edge which is oriented
opposite than towards the root. Notice, however, that this labeling can be
described in an equivalent, simpler way, by saying that all vertices of
$T_{\sigma}$ except from $x_k,x_l,x_m,x_p$ keep their labels in the tree
$T_{\sigma'}$ and the labels of the other vertices can be read by comparing
the Figure \ref{fig:transform0} and Figure \ref{fig:transform2}. 

\begin{figure}
 \psset{unit=2cm}
 \begin{pspicture}[](-1,-0.8)(1,1.2)
 \cnode*(0,0){0.6mm}{A}
 \cnode*(1,0){0.6mm}{C}
 \cnode*(-1,0){0.6mm}{B}
 \ncline[arrowsize=2mm,doubleline=true]{<-}{A}{B}
 \ncline[arrowsize=2mm,doubleline=true]{<-}{A}{C}
 \rput[t](0,-0.1){$x_k$}
 \rput[t](0.9,-0.1){$x_m$}
 \rput[t](-1,-0.1){$x_p$}

\psset{linecolor=gray}

\renewcommand{\psedge}[2]{\ncline[doubleline=true,arrowsize=2mm]{->}{#1}{#2}}
 \rput[B](0,0){\pstree[treemode=U,levelsep=0.5,treesep=0.4]{\Tc*{0.6mm}}{
     \pstree{ \Tc{1.2mm}~[tnpos=l,tnsep=2mm]{$r$} \Tc*{0.6mm}  }{
\renewcommand{\psedge}[2]{\ncline[doubleline=true,arrowsize=2mm]{<-}{#1}{#2}}
   \Tc*{0.6mm}  \Tc*{0.6mm}}
   \renewcommand{\psedge}[2]{\ncline[doubleline=true,arrowsize=2mm]{<-}{#1}{#2}}
   \Tc*{0.6mm}\pstree{\Tc*{0.6mm}}{\Tc*{0.6mm} \Tc*{0.6mm}}}}

  \renewcommand{\psedge}[2]{\ncline[doubleline=true,arrowsize=2mm]{<-}{#1}{#2}}
 \rput[B](1.5,0){\pstree[treemode=R,levelsep=0.5,treesep=0.4]{\Tc*{0.6mm}}{
     \pstree{\Tc*{0.6mm}}{\Tc*{0.6mm}
     \Tc*{0.6mm}} \pstree{\Tc*{0.6mm}}{\Tc*{0.6mm} \Tc*{0.6mm}}}}

 \end{pspicture}
\caption[]{Tree for which the order of the vertices given by $\Delta$
does not coincide with the preorder, the case when $x_l=x_p$.}
 \label{fig:transform0bis}
\end{figure}
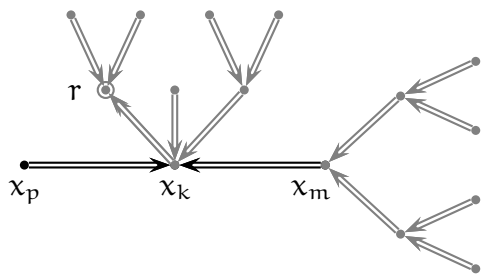

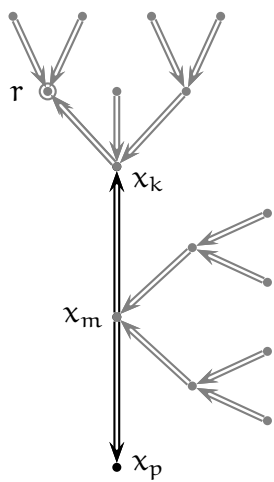
\begin{figure}
\psset{unit=2cm}
\begin{pspicture}[](-1,-1.1)(1,2)
 \cnode*(0,1){0.6mm}{A}
 \cnode*(0,0){0.6mm}{B}
 \cnode*(0,-1){0.6mm}{S}

 \ncline[arrowsize=2mm,doubleline=true]{<-}{A}{B}
 \ncline[arrowsize=2mm,doubleline=true
]{<-}{S}{B}
 \rput[l](0.1,0.9){$x_k$}
 \rput[r](-0.1,0){$x_m$}
 \rput[l](0.1,-1){$x_p$}

\psset{linecolor=gray}
\renewcommand{\psedge}[2]{\ncline[doubleline=true,arrowsize=2mm]{->}{#1}{#2}}
 \rput[B](0,1){\pstree[treemode=U,levelsep=0.5,treesep=0.4]{\Tc*{0.6mm}}{
     \pstree{ \Tc{1.2mm}~[tnpos=l,tnsep=2mm]{$r$} \Tc*{0.6mm}  }{
    
\renewcommand{\psedge}[2]{\ncline[doubleline=true,arrowsize=2mm]{<-}{#1}{#2}}
   \Tc*{0.6mm}  \Tc*{0.6mm}}
   \renewcommand{\psedge}[2]{\ncline[doubleline=true,arrowsize=2mm]{<-}{#1}{#2}}
   \Tc*{0.6mm}\pstree{\Tc*{0.6mm}}{\Tc*{0.6mm} \Tc*{0.6mm}}}}

  \renewcommand{\psedge}[2]{\ncline[doubleline=true,arrowsize=2mm]{<-}{#1}{#2}}
 \rput[B](0.5,0){\pstree[treemode=R,levelsep=0.5,treesep=0.4]{\Tc*{0.6mm}}{
     \pstree{\Tc*{0.6mm}}{\Tc*{0.6mm}
     \Tc*{0.6mm}}  \pstree{\Tc*{0.6mm}}{\Tc*{0.6mm} \Tc*{0.6mm}}}}
 \end{pspicture}
 \caption{%
%
The tree from Figure \ref{fig:transform0bis} after regluing the edges
$(x_l,x_k)$ and $(x_m,x_k)$ in a different way. Notice the change of the
labels of
the vertices $x_k,x_m,x_p$.}
 \label{fig:transform2bis}
\end{figure}

For simplicity, in the above discussion we considered only the case when
$x_l\neq x_p$; otherwise only a minor correction is necessary, namely Figure
\ref{fig:transform0} must be replaced by Figure \ref{fig:transform0bis} and
Figure \ref{fig:transform2} by
Figure \ref{fig:transform2bis}.

In this way our analysis of case (\ref{case:bad-case-1}) is finished. Notice
that in this case
$\big|\{v \in T_{\sigma}: v \succ r\}\big| > 
\big|\{v\in T_{\sigma'}: v \succ r\}\big|$ therefore in the graph $\graf$
there is an oriented edge pointing from the vertex correspondning to the
equivalence class $(T_{\sigma},[\Delta]_{\sim_d})$ to the vertex corresponding
to the equivalence class $(T_{\sigma'},[\Delta']_{\sim d})$.

In the case (\ref{case:bad-case-2}) involution $f$ is just the inverse of the
map $f$ described above for the case (\ref{case:bad-case-1}); in this case the
edge in the graph $\graf$ is oriented in the opposite direction as in the case
(\ref{case:bad-case-1}).


\begin{figure}
\psset{unit=2cm}
\begin{pspicture}[](-1,-1.1)(1,2)
 \cnode*(0,1){0.6mm}{A}
 \cnode*(0,0){0.6mm}{B}
 \cnode*(0,-1){0.6mm}{S}

 \ncline[arrowsize=2mm,doubleline=true]{<-}{A}{B}
 \ncline[arrowsize=2mm,doubleline=true
]{<-}{S}{B}
 \rput[l](0.1,0.9){$x_k$}
 \rput[l](0.2,0){$x_l$}
 \rput[l](0.1,-1){$x_p$}

\psset{linecolor=gray}
\renewcommand{\psedge}[2]{\ncline[doubleline=true,arrowsize=2mm]{->}{#1}{#2}}
 \rput[B](0,1){\pstree[treemode=U,levelsep=0.5,treesep=0.4]{\Tc*{0.6mm}}{
     \pstree{ \Tc{1.2mm}~[tnpos=l,tnsep=2mm]{$r$} \Tc*{0.6mm}  }{
    
\renewcommand{\psedge}[2]{\ncline[doubleline=true,arrowsize=2mm]{<-}{#1}{#2}}
   \Tc*{0.6mm}  \Tc*{0.6mm}}
   \renewcommand{\psedge}[2]{\ncline[doubleline=true,arrowsize=2mm]{<-}{#1}{#2}}
   \Tc*{0.6mm}\pstree{\Tc*{0.6mm}}{\Tc*{0.6mm} \Tc*{0.6mm}}}}

 \renewcommand{\psedge}[2]{\ncline[doubleline=true,arrowsize=2mm]{<-}{#1}{#2}}
 \rput[B](-0.55,0){\pstree[treemode=L,levelsep=0.5,treesep=0.4]{\Tc*{0.6mm}}{
     \pstree{\Tc*{0.6mm}}{\Tc*{0.6mm} \Tc*{0.6mm}}
\Tc*{0.6mm}\pstree{\Tc*{0.6mm}}{\Tc*{0.6mm}
     \Tc*{0.6mm}~[tnpos=l,tnsep=0mm]{$x_m$}    }}}

 \rput[B](0.5,0){\pstree[treemode=R,levelsep=0.5,treesep=0.4]{\Tc*{0.6mm}}{
     \pstree{\Tc*{0.6mm}}{\Tc*{0.6mm}
     \Tc*{0.6mm}}  \pstree{\Tc*{0.6mm}}{\Tc*{0.6mm} \Tc*{0.6mm}}}}
 \end{pspicture}
 \caption{The tree from Figure \ref{fig:transform0} after regluing the edges
$(x_l,x_k)$ and $(x_m,x_k)$ in a different way. Notice the change of the
labels of
the vertices $x_k,x_l,x_m,x_p$. }

 \label{fig:transform2}
\end{figure}
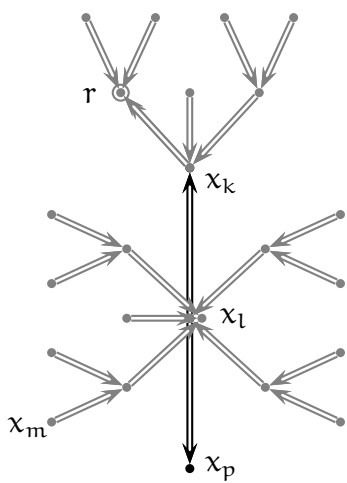


\subsection{Proof of Proposition \ref{prop:sinks_exist}}
\begin{proof}[Proof of Proposition \ref{prop:sinks_exist}]
Any configurations $(T_{\sigma},\Delta_1)$,
$(T_{\sigma},\Delta_2)$ which contribute to the same vertex of $\graf$ can be
transformed into each other by repeatedly interchanging the order of adjacent
derivatives which remove an edge oriented towards the root and an edge oriented
opposite to towards the root. If these configurations give rise to outgoing
edges in the graph $\graf$ it follows that the case (\ref{case:bad-case-1})
holds true for both $(T_{\sigma},\Delta_1)$ and $(T_{\sigma},\Delta_2)$. One
can easily see that in both cases the procedure described in Section
\ref{subsec:explicit_form} unglues and reglues the same two edges hence the
resulting configurations $f(T_{\sigma},\Delta_1)$ and $f(T_{\sigma},\Delta_2)$
belong to the same equivalence class hence correspond to the same vertex of
$\graf$. This shows that every vertex of $\graf$ has at most one outgoing edge.

If the vertices corresponding to the configurations $(T_{\sigma},\Delta)$,
$(T_{\sigma'},\Delta')$ are connected by an oriented edge then 
$$ \big| \{ x\in T_{\sigma}: x \succ r \} \big| > \big| \{ x\in T_{\sigma'}: x
\succ r \} \big|$$ 
therefore there are no oriented cycles in $\graf$.

The second part of Proposition \ref{prop:sinks_exist} follows easily from the
description of map $f$ in Section \ref{subsec:explicit_form}.
\end{proof}

\section{Acknowledgments} 
This research work was supported by the Ministry of Education and Science in Poland in years 2006-2009 under grant number 1 P03A 013 30. Research was partially supported by EU Research Training Network `QP-Applications', contract HPRN-CT-2002-00279. Research was partially supported by European Commission Marie Curie Host Fellowship for the Transfer of Knowledge `Harmonic Analysis, Nonlinear Analysis and Probability' grant MTKD-CT-2004-013389.

\bibliographystyle{alpha}
\bibliography{biblio}
\end{document}